\documentclass[10pt,leqno]{article}

\usepackage{amssymb} 
\usepackage{amsmath} 

\input xy 
\xyoption{all}

\title{Classification of linearly compact simple Jordan
and
generalized Poisson superalgebras} 
\author{{\sc Nicoletta Cantarini}\thanks{Dipartimento di Matematica
Pura ed Applicata, Universit\`a di Padova, Padova, Italy}
\and\setcounter{footnote}{6}
{\sc Victor G.\ Kac}\thanks{Department of Mathematics, MIT, Cambridge,
Massachusetts 02139, USA, and IHES, 91440 Bures-sur-Yvette, France}}
\date{\small{\em To Ernest Borisovich Vinberg on his 70th birthday}}

\newtheorem{theorem}{Theorem}[section] 
\newtheorem{lemma}[theorem]{Lemma} 
\newtheorem{corollary}[theorem]{Corollary} 
\newtheorem{proposition}[theorem]{Proposition} 
\newtheorem{definition}[theorem]{Definition} 
\newtheorem{remark}[theorem]{Remark}
\newtheorem{remarks}[theorem]{Remarks}
\newtheorem{example}[theorem]{Example}
\newtheorem{examples}[theorem]{Examples}
\newtheorem{question}[theorem]{Question}
\def\Z{\mathbb{Z}} 

\def\g{\mathfrak{g}}
\def\a{\mathfrak{a}}
\def\b{\mathfrak{b}}

\def\h{\mathfrak{h}} 
\def\t{\mathfrak{t}} 

\def\C{\mathbb{C}}
\def\F{\mathbb{F}}

\def\H{\mathbb{H}}

\numberwithin{equation}{section}

\makeatletter

\def\enumerate{%
  \ifnum \@enumdepth >\thr@@\@toodeep\else
    \advance\@enumdepth\@ne
    \edef\@enumctr{enum\romannumeral\the\@enumdepth}%
      \list
        {\csname label\@enumctr\endcsname}%
        {\usecounter\@enumctr
          \addtolength{\leftmargin}{-\leftmargin}
          \settowidth{\labelwidth}{(99)}
          \itemindent = \labelwidth
          \addtolength{\itemindent}{\labelsep}
        \listparindent=1em      
          \def\makelabel##1{{##1}\hfill}
          }%
  \fi}

\makeatother

\begin{document} 
\maketitle 
\date
\begin{abstract} 
We classify all linearly compact simple Jordan
superalgebras over an algebraically closed
field of characteristic zero. As a corollary, we deduce the classification of 
all linearly compact unital simple generalized Poisson superalgebras.
\end{abstract} 
\section*{Introduction}
Jordan algebras first appeared in a 1933 paper by P.\ Jordan on
axiomatic foundations of quantum mechanics. Subsequently,
Jordan, von Neumann and Wigner classified the ``totally positive''
simple finite-dimensional Jordan algebras \cite{JNW}.
All simple finite-dimensional Jordan algebras over an algebraically
closed field $\mathbb{F}$ of characteristic different from 2 have been
classified some years later by Albert \cite{A}, using the idempotent method.

Albert's list is as follows:
\begin{itemize}
\item[A.] $gl(m)_+$: the space of $m\times m$-matrices over $\mathbb{F}$
with multiplication $a\circ b=\frac{1}{2}(ab+ba)$;
\item[B.] $o(m)_+$: the subalgebra of symmetric matrices in $gl(m)_+$;
\item[C.] $sp(2m)_+$: the subalgebra of matrices in $gl(2m)_+$ fixed 
under the anti-involution $a\mapsto J^{-1}a^tJ$, where $J$ is a non-degenerate
skewsymmetric matrix;
\item[D.] $(m)_+$: $\mathbb{F}e\oplus\mathbb{F}^m$, where $e$ is the
unit element and the multiplication on $\mathbb{F}^m$ is defined
via a non-degenerate symmetric bilinear form $(\cdot,\cdot)$:
$a\circ b=(a,b)e$;
\item[E.] the 27-dimensional Jordan algebra $E$ of Hermitian
$3\times 3$-matrices over Cayley numbers.
\end{itemize}

Jordan, von Neumann and Wigner expressed a hope in their above
cited paper that there should be new types of simple infinite-dimensional
Jordan algebras, but in his remarkable paper \cite{Z} Zelmanov proved
that this was not the case. Namely, all examples $A_+$ of type $A$
can be obtained by replacing $m\times m$-matrices by an infinite-dimensional
simple associative algebra $A$; all examples of types $B$ and $C$
can be obtained by taking the subalgebra of the Jordan algebra $A_+$,
fixed under an anti-involution of $A$; all examples of type $D$ are obtained
by replacing $\mathbb{F}^m$ by an infinite-dimensional vector space
with a non-degenerate symmetric bilinear form;
as a result, along with $E$, one gets a complete list of simple
Jordan algebras over $\F$. 

One of the advances of the theory of Jordan algebras in the 60s
was the discovery of the Tits-Kantor-Koecher (TKK) construction
\cite{T}, \cite{Ka1}, \cite{Ko}. It is based on the observation
that if $\g=\g_{-1}\oplus\g_0\oplus\g_1$ is a Lie algebra with 
a {\em short} $\Z$-grading and $f\in\g_1$, then the formula
\begin{equation}
a\circ b=[[f,a],b]
\label{basicidea}
\end{equation}
defines a structure of a Jordan algebra on $\g_{-1}$. If $e\in\g_{-1}$
is the unit element of this Jordan algebra and $h=[e,f]$ for some
$f\in\g_1$, then
$\{e, h, f\}$ is an $sl_2$-triple in $\g$, and the eigenspace
decomposition of $ad \,h$ defines a short grading of $\g$. This leads
to an equivalence of the category of unital Jordan algebras
to the category of pairs $(\g,\a)$, where $\g$ is a Lie algebra
and $\a$ is its $sl_2$ subalgebra whose semisimple element with
eigenvalues $0,1,-1$ defines a short grading of $\g$, such that
$\g_0$ contains no non-zero ideals of $\g$. (Such $\a$ is called
a short subalgebra of $\g$.) 
The ($\Z$-graded) Lie
algebra corresponding to a Jordan algebra $J$ under this
equivalence is called the TKK construction of $J$.
This  construction can be applied to non-unital
Jordan algebras as well, but then the corresponding category
of $\Z$-graded Lie algebras is a bit more complicated \cite{K1}. However,
it is easier to add the unit element to the Jordan algebra if
it does not have one, which leads to the $sl_2$ algebra of
outer derivations of the corresponding Lie algebra. In the present paper
we adopt the latter approach.

Using the TKK construction, Kac \cite{K1} derived a classification
of simple finite-dimensional Jordan superalgebras, under the assumption
char $\mathbb{F}=0$, from the classification of simple finite-dimensional
Lie superalgebras \cite{K2}. The list consists of Jordan superalgebras of
types $A$, $BC$ and $D$, which are super analogues of the types 
$A$, $B$, $C$, $D$
above, two ``strange'' series $P$ and $Q$, the Kantor series
\cite{Ka2}, and, in addition to the exceptional Jordan algebra $E$,
a family of 4-dimensional Jordan superalgebras $D_t$, a non-unital
3-dimensional Jordan superalgebra $K$, and an exceptional
10-dimensional Jordan superalgebra $F$. (Note that
the multiplication table of $F$ in \cite{K1} contained errors,
which were fixed in \cite{HK}, and that the Kantor series 
was inadvertently omitted
in \cite{K1}, as it has been pointed out in \cite{Ka2}.)

It turned out that the Kantor series is a subseries of
finite-dimensional Jordan superalgebras of a family
of simple Jordan superalgebras, obtained
via the Kantor-King-McCrimmon (KKM) double, defined as follows.
Let $A$ be a commutative associative superalgebra, let $J(A)=A\oplus \eta A$,
where $\eta$ is an odd indeterminate, and define a product on
$J(A)$ by $(a, b\in A)$:
\begin{equation}
a\circ b=ab, ~~\eta a\circ b=\eta(ab), ~~ a\circ \eta b=(-1)^{p(a)}\eta(ab), 
~~\eta a\circ \eta b=(-1)^{p(a)}\{a,b\},
\label{0.2}
\end{equation}
where $\{\cdot ,\cdot\}$ is some other bilinear product on $A$.
Kantor \cite{Ka2} checked that if $\{\cdot, \cdot\}$ is a Poisson
bracket, i.e., a Lie superalgebra bracket satisfying the Leibniz rule,
then $J(A)$ is a Jordan superalgebra. For example, the
Grassmann superalgebra $\Lambda(n)=\Lambda(\xi_1, \dots, \xi_n)$
with bracket \cite{K2}:
\begin{equation}
\{f,g\}=(-1)^{p(f)}\sum_{i=1}^n\frac{\partial f}{\partial\xi_i}
\frac{\partial g}{\partial\xi_i}
\label{bracketH(n)}
\end{equation}
is a Poisson superalgebra. Hence $JP(0,n)=\Lambda(n)\oplus\eta\Lambda(n)$
with product (\ref{0.2}) is a Jordan superalgebra, and this is
the Kantor series, missed in \cite{K1}.

King and Mc Crimmon \cite{KMCC} found necessary and sufficient conditions
on the bracket $\{\cdot, \cdot\}$ in order for (\ref{0.2}) to give a
Jordan product. Using this, we show that all the KKM doubles $J(A)$, which are
Jordan superalgebras, can be obtained from a
generalized Poisson bracket on the superalgebra $A$. By definition, this is a
Lie superalgebra
bracket $\{\cdot,\cdot\}$, satisfying the generalized Leibniz rule:
\begin{equation}
\{a,bc\}=\{a,b\}c+(-1)^{p(a)p(b)}b\{a,c\}+D(a)bc,
\label{genLeibniz}
\end{equation}
where $D$ is an even derivation of the product and the bracket. 
(If $D=0$ we get the
usual Leibniz rule, and if $A$ has a unit element $e$, then $D(a)=\{e, a\}$.)

Now, the commutative associative superalgebra ${\cal O}(m,n)=\Lambda(n)[[x_1,
\dots, x_m]]$ carries a well known structure of a generalized Poisson 
superalgebra \cite{K2}, \cite{K3}, which we denote by $P(m,n)$.
The KKM double of this is a simple (unless $m=n=0$) 
Jordan superalgebra, denoted by 
$JP(m,n)$. The case when $m=1$ and ${\cal O}(1,n)$ is replaced by
$\Lambda(n)[x,x^{-1}]$ plays a prominent role in the Kac-Martinez-Zelmanov
(KMZ) classification \cite{KMZ}; denote this Jordan superalgebra
by $JP(1,n)^{\sim}$. It is shown in \cite{KMZ}, that $JP(1,3)^{\sim}$
contains a ``half size'' exceptional Jordan superalgebra $JCK^{\sim}$,
whose TKK construction is the exceptional Lie superalgebra $CK_6$, discovered by
Cheng and Kac \cite{CK0}. The KMZ classification theorem
\cite{KMZ} states that a complete list of $\Z$-graded simple
unital Jordan superalgebras of growth 1 consists of twisted
loop algebras over finite-dimensional simple Jordan superalgebras,
Jordan superalgebras of a non-degenerate supersymmetric bilinear form
on a $\Z$-graded vector superspace (type $D$ above), the series 
$JP(1,n)^{\sim}$, $JCK^{\sim}$ and their twisted analogues, and 
simple Jordan superalgebras of ``Cartan type''.

By definition, a Lie superalgebra $L$ is called of Cartan type if it has
a subalgebra $L_0$ of finite codimension,
containing no non-zero ideals of $L$; a Jordan superalgebra is of
Cartan type if it contains a subalgebra of finite codimension and
its TKK construction is of Cartan type. One can construct a
descending filtration of a Lie superalgebra of Cartan type
$L\supset L_0\supset L_1\supset L_2\supset\dots$ (see e.g.\ \cite{K3}),
introduce  a structure of a topological Lie superalgebra on $L$
by taking $\{L_j\}_{j\in\Z_+}$ to be the fundamental system
of neighborhoods of zero, and consider the completion $\overline{L}$
of $L$ with respect to this topology. The topological Lie superalgebra
$\overline{L}$ is linearly compact in the sense of the following definition.

A topological superalgebra (associative, or Lie, or Jordan) is called
linearly compact if the underlying topological space is linearly
compact (i.e., it is isomorphic to the topological direct product
of finite-dimensional vector spaces with the discrete topology). It is
well known (see \cite{G1})
that any linearly compact Lie superalgebra has a subalgebra
of finite codimension, hence any linearly compact Lie superalgebra
is of Cartan type.

Infinite-dimensional simple linearly compact Lie superalgebras
were classified by Kac \cite{K3}. The main goal of the present
paper is to establish the equivalence (similar to the discussed above)
of the category of unital linearly compact Jordan superalgebras
and the category of linearly compact Lie superalgebras with a short
subalgebra $\a=\langle e,h,f\rangle$ such that the centralizer of
$h$ contains no non-zero ideals of the whole algebra, and to use this
equivalence in order to derive the classification of infinite-dimensional
linearly  compact simple Jordan superalgebras from the 
classification in the Lie superalgebra case.

It turns out that among the ten series of infinite-dimensional linearly
compact simple Lie superalgebras only the series $H(m,n)$ and
$K(m,n)$ admit a (unique) short subalgebra (provided that $n\geq 3$),
and the superalgebra $S(1,2)$ admits a short algebra of outer derivations.
In addition, among the five exceptional ones, only $E(1,6)$ admits a
(unique) short subalgebra.
As a result, we obtain the following list of all infinite-dimensional
linearly compact simple Jordan superalgebras:
 the series of unital Jordan superalgebras
$JP(m,n)$ (mentioned above), the exceptional unital Jordan
superalgebra $JCK$ (which is a ``half-size'' subalgebra of
$JP(1,3)$), and a new, non-unital Jordan superalgebra, which we denote by
$JS$.

The construction of the latter is very simple. Let $JS={\cal O}(1,1)$
with reversed parity, and the product:
\begin{equation}
a\circ b=aD(b)+(-1)^{p(a)}D(a)b,
\label{newforintro}
\end{equation}
where $D=\frac{\partial}{\partial\xi}+\xi\frac{\partial}{\partial x}.$

Along the way, we redo the (corrected) classification of \cite{K1}, using the 
simpler method of short subalgebras.

As a corollary, we obtain that any simple unital generalized Poisson
superalgebra is isomorphic to one of the 
$P(m,n)^{\varphi}$, 
where $P(m,n)$ is one of the ``standard'' generalized Poisson superalgebras,
mentioned above, and 
$P(m,n)^{\varphi}$ is obtained from it by twisting the bracket by an even 
invertible element $\varphi \in {\cal O}(m,n)$ by the formula:
\begin{equation}
\{a,b\}^{\varphi}=\varphi^{-1}\{\varphi a,\varphi b\}.
\end{equation}

We hope that the blow to Jordan-von Neumann-Wigner expectations, dealt by
Zelmanov's theorem on non-existence of new types of simple
infinite-di\-mensional Jordan algebras, will be remedied by applications
of simple infinite-dimensional Jordan superalgebras to supersymmetric quantum
mechanics.

We are grateful to S.-J.\ Cheng and A.\ A.\ Kirillov for correspondence.

Throughout the paper the base field $\F$ is algebraically closed
of characteristic 0, unless otherwise specified.

\section{Preliminaries on superalgebra}
We recall that a vector {\em superspace} $V$ is a $\Z/2\Z$-graded vector space:
$V=V_{\bar{0}}\oplus V_{\bar{1}}$. 
If $a$ lies in $V_{\alpha}$, $\alpha\in\Z/2\Z=\{\bar{0}, \bar{1}\}$, we say that
$\alpha$ is the {\em parity} of $a$ and we denote it by $p(a)$.
The elements of $V_{\bar{0}}$ are called
{\em even}, the elements of $V_{\bar{1}}$ are called {\em odd}.

A bilinear form $(\cdot, \cdot)$ on a vector superspace $V$ is called {\em supersymmetric}
(resp.\ {\em superskew-symmetric}) if its restriction to
$V_{\bar{0}}$ is symmetric (resp.\ skew-symmetric), its
restriction to $V_{\bar{1}}$ is skew-symmetric (resp.\ symmetric)
and \break $(V_{\bar{0}}, V_{\bar{1}})=0$. The matrix of  a supersymmetric
(resp.\ superskew-symmetric) bilinear form on a vector superspace $V$ 
in some
basis of $V$,
is called supersymmetric (resp.\ superskew-symmetric).

A {\em superalgebra} is a $\Z/2\Z$-graded algebra.
A superalgebra $A$ is called {\em commutative} if, for every $a,b\in A$,
$$ab=(-1)^{p(a)p(b)}ba.$$
Associativity of superalgebras is defined as for algebras.

Given a vector superspace $V$, the associative algebra $End(V)$ is
equipped with the induced $\Z/2\Z$-grading: $End(V)=\oplus_{\alpha\in\Z/2\Z}
End_{\alpha}(V)$, where
$$End_{\alpha}(V)=\{f\in End(V)~|~f(V_s)\subseteq V_{s+\alpha}\}.$$
The commutator in an associative superalgebra is defined by
$$[a,b]=ab-(-1)^{p(a)p(b)}ba.$$
With respect to this bracket, $End(V)$ is a Lie superalgebra.

A derivation with parity $s$, $s\in\Z/2\Z$, of a superalgebra $A$
is an endomorphism $D\in End_s(A)$ with the property:
$$D(ab)=D(a)b+(-1)^{p(a)s}aD(b).$$
The space $Der A$ of all derivations of $A$ is a subalgebra of
$End(A)$.

A super anti-involution of a superalgebra $A$ is a grading preserving
linear map $^*: A\rightarrow A$ such that
$(a^*)^*=a$ for every $a\in A$, and
$$(ab)^*=(-1)^{p(a)p(b)}b^*a^* ~~(a,b\in A).$$
If $(\cdot, \cdot)$ is a non-degenerate supersymmetric bilinear form on a 
vector superspace
$V$ of dimension $(m|n)$, then we have a super anti-involution 
$^*: End(V)\rightarrow End(V)$,
defined by
\begin{equation}
(f(a),b)=(-1)^{p(a)p(f)}(a, f^*(b)),~~f\in End(V), ~a,b \in V.
\label{anti-inv0}
\end{equation}
If $B$ is the matrix of this bilinear form
in some basis of $V$, compatible with grading, 
and $M$ (resp.\ $M^*$) is the matrix of an endomorphism $f$ (resp.\ $f^*$)
in the same basis, then
\begin{equation}
a^tM^tBb=(-1)^{p(a)p(M)}a^tBM^*b,
\label{anti-inv}
\end{equation}
where $^*: Mat(m,n) \rightarrow Mat(m,n)$
is a super anti-involution given by: $M^*=B^{-1}M^TB$, and 
$$\left(\begin{array}{cc}
\alpha & \beta\\ 
\gamma & \delta 
\end{array} 
\right)^T=\left(\begin{array}{cc} 
\alpha^t & \gamma^t\\ 
-\beta^t & \delta^t 
\end{array} 
\right),$$ 
with $^t$ standing for the usual transpose.

If $(\cdot, \cdot)$ is an odd supersymmetric bilinear form, i.e.,
its matrix in some basis of $V$, compatible with grading, is
$\left(\begin{array}{cc}
0& I\\
I& 0
\end{array}
\right)$, then the super anti-involution $^*$, defined
in (\ref {anti-inv0}), in matrix form  looks as (\ref{anti-inv}),
where $^*: Mat(m,m) \rightarrow Mat(m,m)$
is given by:
$$\left(\begin{array}{cc}
\alpha & \beta\\
\gamma & \delta
\end{array}
\right)^*=
\left(\begin{array}{cc}
\delta^t & \beta^t\\
-\gamma^t & \alpha^t
\end{array}
\right).$$

Given a vector superspace $V$, 
a bilinear map $\varphi: V\times V\longrightarrow
V$ is called commutative if $\varphi(a,b)=(-1)^{p(a)p(b)}\varphi(b,a)$;
it is called anticommutative if 
$\varphi(a,b)=-(-1)^{p(a)p(b)}\varphi(b,a)$.

\section{Preliminaries on linearly compact spaces and algebras}
Recall that a {\em linear topology} on a vector space $V$ over $\F$
is a topology for which there exists a fundamental system of neighbourhoods
of zero, consisting of vector subspaces of $V$. A vector space
$V$ with linear topology is called {\em linearly compact} if there
exists a fundamental system of neighbourhoods of zero, consisting of subspaces
of finite codimension in $V$ and, in addition, $V$ is complete in this
topology (equivalently, if $V$ is a topological direct product of some number
of copies of $\F$ with discrete topology).

If $\varphi: V \rightarrow U$ is a continuous map
of vector spaces with linear topology and $V$ is linearly
compact, then $\varphi(V)$ is linearly compact \cite{G1}.

The basic examples of linearly compact spaces are finite-dimensional
vector spaces with the discrete topology, and the space of formal
power series $V[[x_1, \dots, x_m]]$ over a finite-dimensional vector 
space $V$, with
the topology, defined by taking as a fundamental system of
neighborhoods of $0$ the subspaces $\{x_1^{j_1}\dots x_m^{j_m}
V[[x_1, \dots, x_m]]\}_{(j_1, \dots, j_m)\in\Z_+^m}$. 

Let $U$ and $V$ be two vector spaces with linear topology, and let
$Hom(U,V)$ denote the space of all continuous linear maps from $U$ to $V$.
We endow the vector space $Hom(U,$ $V)$ with a ``compact-open'' linear 
topology, the fundamental system of neighbourhoods of zero
being $\{\Omega_{K,W}\}$, where $K$ runs over all linearly compact
subspaces of $U$, $W$ runs over all open subspaces of $V$,
and 
$$\Omega_{K,W}=\{\varphi\in Hom(U,V)~|~\varphi(K)\subset W\}.$$
In particular, we define the
dual space of $V$ as $V^*=Hom(V,\F)$, where $\F$ is endowed with the 
discrete topology. It is easy to see that $V$ is linearly compact
if and only if $V^*$ is discrete. Note also that a discrete space $V$
is linearly compact if and only if $\dim V<\infty$.

The tensor product of two vector spaces $U$ and $V$ with linear topology
is defined as
$$U\hat{\otimes} V=Hom(U^*, V).$$
Thus, $U\hat{\otimes} V=U\otimes V$ if both $U$ and $V$ are discrete
and $U\hat{\otimes} V=(U^*\otimes V^*)^*$ if both $U$ and $V$ are
linearly compact. Hence the tensor product of linearly compact spaces
is linearly compact \cite{G1}.

Recall that the basic property of the ``compact-open'' topology on $Hom(U,V)$
is the following canonical isomorphism of vector spaces with linear
topology:
\begin{equation}
Hom(M, Hom(U,V))\simeq Hom(M\hat{\otimes} U, V).
\label{compactopentop}
\end{equation}
(This isomorphism is the closure of the map $\Phi(\psi)(m\otimes u)=
\psi(m)(u)$.)

A {\it linearly
compact superalgebra} is a topological superalgebra whose underlying
topological space is linearly compact. Of course, any
finite-dimensional superalgebra is linearly compact. The basic example of an
associative linearly compact superalgebra is the 
superalgebra
${\cal O}(m,n)=\Lambda(n)[[x_1,\dots, x_m]]$, where $\Lambda(n)$
denotes the Grassmann algebra on $n$ anticommuting indeterminates
$\xi_1, \dots, \xi_n$, and the superalgebra parity is defined by
$p(x_i)=\bar{0}$, $p(\xi_j)=\bar{1}$, with the formal topology defined
as above for $V=\Lambda(n)$.
The basic example of a
linearly compact Lie superalgebra is $W(m,n)=Der {\cal O}(m,n)$, the
Lie superalgebra of all continuous derivations of the superalgebra
${\cal O}(m,n)$. One has:
$$W(m,n):=\{X=\sum_{i=1}^mP_i(x,\xi)\frac{\partial}{\partial
x_i}+\sum_{j=1}^nQ_j(x,\xi)\frac{\partial}{\partial\xi_j} ~|~ P_i,
Q_j\in{\cal O}(m,n)\}.$$
\begin{lemma}\label{lincomp} Let $A$ be a linearly compact (super)algebra, and let
$L_a$ denote the operator of left multiplication by $a$. Then:
\begin{enumerate}
\item[$(a)$] the subspace $\{L_a~|~a\in A\}$ of the space of all 
continuous linear maps $Hom(A,A)$ is linearly compact;
\item[$(b)$] the closure of the set $[L_a, L_b]$ in $Hom(A,A)$ is linearly compact.
\end{enumerate}
\end{lemma}
{\bf Proof.} We apply (\ref{compactopentop}) to $M=A$ (resp.\
$M=A\hat{\otimes} A$) and $U=V=A$ for the proof of $(a)$ (resp.\ $(b)$).
It follows that the map $a\mapsto L_a$ is continuous, and since the image of
a linearly compact space under a continuous map is linearly compact,
 $(a)$ follows. The proof of $(b)$ is similar by making use of the
fact that $A\hat{\otimes} A$ is linearly compact. \hfill$\Box$

\section{Generalized Poisson superalgebras and the\break Schouten bracket}
\begin{definition}
A generalized Poisson superalgebra is a superalgebra $A$ with 
two operations: 
a commutative associative superalgebra product $(a,b)\mapsto ab$ and a Lie
superalgebra bracket $(a,b)\mapsto \{a,b\}$, jointly satisfying the
generalized Leibniz
rule, namely, for $a, b, c\in A$, one has
\begin{equation}
\{a, bc\}=\{a,b\}c+(-1)^{p(a)p(b)}b\{a,c\}+D(a)bc,
\label{gP}
\end{equation} 
for some even derivation $D$ of $A$ with respect to the product and the 
bracket.
If $D=0$, then relation (\ref{gP}) becomes the usual Leibniz rule; in this case
$A$ is called a Poisson superalgebra. A generalized Poisson superalgebra is 
called simple if it contains no non-trivial ideals with respect to both 
operations and the bracket is non-trivial. 
\end{definition}

\begin{remark}\label{unitalgP}\em If $A$ is a unital generalized Poisson 
superalgebra 
and $e$ is its unit element, then $D(a)=\{e,a\}$. Indeed, it is sufficient 
to set
$b=c=e$ in (\ref{gP}).
\end{remark}

\begin{example}\label{P(2k,n)}\em
Consider the associative superalgebra ${\cal O}(2k,n)$
with even indeterminates $p_1, \dots, p_k, q_1,\dots,
q_k$, and odd indeterminates $\xi_1,\dots, \xi_n$.
Define on ${\cal O}(2k,n)$ the following bracket
($f,g\in {\cal O}(2k,n)$):
\begin{equation}\{f,g\}=\sum_{i=1}^k(\frac{\partial f}{\partial p_i}\frac{\partial
g}{\partial q_i}-\frac{\partial f}{\partial q_i}\frac{\partial
g}{\partial p_i})+(-1)^{p(f)}(\sum_{i=1}^{n-2}\frac{\partial f}{\partial \xi_i}
\frac{\partial g}{\partial \xi_i}+
\frac{\partial f}{\partial \xi_{n-1}}
\frac{\partial g}{\partial \xi_n}+
\frac{\partial f}{\partial \xi_n}
\frac{\partial g}{\partial \xi_{n-1}}).
\label{bracket}
\end{equation}
Then ${\cal O}(2k,n)$ with this bracket is a Poisson superalgebra
(cf.\ \cite[\S 1.2]{CK}), which we
denote  by $P(2k,n)$.

Consider ${\cal O}(2k,n)$ with its Lie superalgebra structure defined by
bracket (\ref{bracket}). Then
$\F 1$ is a central ideal, hence
$H'(2k,n):={\cal O}(2k,n)/\F 1$ is a linearly compact Lie 
superalgebra \cite{CK}. Set 
$H(2k,n)=[H'(2k,n), H'(2k,n)]$. 
We recall that
$H(2k,n)=H'(2k,n)$ if $k\geq 1$, and then it is a simple
Lie superalgebra;
besides, $H(0,n)$ is simple if and only
if $n\geq 4$ and $H'(0,n)=H(0,n)+\F\xi_1\dots\xi_n$
(see \cite[Proposition 3.3.6]{K1}, \cite[Example 4.3]{K2}).
\end{example}

\begin{example}\label{P(2k+1,n)}\em
Consider the associative superalgebra ${\cal O}(2k+1,n)$
with even indeterminates
$t, p_1, \dots, p_k, q_1,\dots,
q_k$, and odd indeterminates $\xi_1,\dots, \xi_n$.
Define on ${\cal O}(2k+1,n)$ the following bracket
($f,g\in {\cal O}(2k+1,n)$):
\begin{equation}
\{f,g\}_K=(2-E)f\frac{\partial g}{\partial t}-\frac{\partial
f}{\partial t}(2-E)g+\{f,g\},
\label{Leibniz!!}
\end{equation}
where $E=\sum_{i=1}^k(p_i\frac{\partial}{\partial
p_i}+q_i\frac{\partial}{\partial q_i})+\sum_{i=1}^n\xi_i
\frac{\partial}{\partial \xi_i}$ is the Euler operator
and $\{f,g\}$ is given by (\ref{bracket}).
Then ${\cal O}(2k+1,n)$ with bracket $\{\cdot,\cdot\}_K$
 is a generalized Poisson superalgebra
with $D=2\frac{\partial}{\partial t}$ (cf.\ \cite[Remark 2.19]{CantaK}),
which we denote  by $P(2k+1,n)$.  

We shall denote by $K(2k+1,n)$ the superalgebra
$P(2k+1,n)$ with its Lie superalgebra structure.
This is a simple linearly compact Lie superalgebra for every $k$ and $n$
(cf.\ \cite[\S 1.2]{CK}).
\end{example}

\begin{remark}\label{anotherway}\em
Let $C=(f_{ij})$ be a $(2k+n)\times (2k+n)$ superskew-symmetric matrix over 
$\F$.  
Then a Poisson bracket 
can be defined on ${\cal O}(2k,n)$ as follows: denote by $X_i$,
with $i=1,\dots,2k$ and $i=2k+1,\dots, 2k+n$, all even and odd indeterminates,
respectively. Define
\begin{equation}
\{X_i, X_j\}=f_{ij},
\label{simpleway}
\end{equation}
and then extend $\{\cdot,\cdot\}$ to all polynomials by the usual Leibniz
rule. Next, consider the associative superalgebra ${\cal O}(2k+1,n)$
with the same indeterminates as above
and one more even indeterminate $t$. Then a
generalized Poisson bracket, with $D=2\frac{\partial}{\partial t}$,
 can be defined
on ${\cal O}(2k+1,n)$  by (\ref{simpleway}) and 
\begin{equation}
\{t, X_i\}=-X_i,
\label{tcomm}
\end{equation}
extending by the generalized Leibniz rule (\ref{gP}).
Note that Examples \ref{P(2k,n)} and  \ref{P(2k+1,n)} are isomorphic to 
these with a non-degenerate matrix $C$, by an invertible linear change of 
indeterminates. Note also that all generalized Poisson superalgebras $P(m,n)$,
except for $P(0,0)$, are simple. 
\end{remark}

\medskip



Now we briefly discuss the connection of generalized Poisson
superalgebras with the Schouten bracket.

Let $A$ be a commutative associative superalgebra over $\F$.
Then $Der A$ is a Lie superalgebra over $\F$ and a left
$A$-module. Define on
$\Lambda(Der A)$, the exterior algebra over $Der A$,
 viewed as a module over $A$,  the following parity:
$p(fX_1\wedge\dots\wedge X_k)=p(f)+p(X_1)+\dots+p(X_k)+k$,
for $f\in A$ and $X_i\in Der A$, $i=1, \dots,k$. 
The Schouten bracket $[\cdot, \cdot]_S$ on
$\Lambda(Der A)$ is
defined as follows:
$$[f,g]_S=0 ~~\mbox{for}~ f,g\in A; ~~~[X,Y]_S=[X,Y] ~~\mbox{for}~ X, Y\in Der A;$$
$$[D,f]_S=-(-1)^{(p(D)+1)(p(f)+1)}[f,D]_S=D(f) ~~\mbox{for}~ D\in Der A, f\in A;$$
and extended by the odd Leibniz rule:
$$[a,b\wedge c]_S=[a,b]_S\wedge c+(-1)^{(p(a)+1)p(b)}b\wedge[a,c]_S.$$
Then $[\cdot,\cdot]_S$ satisfies the Lie superalgebra axioms after
reversing parity.
\begin{example}\em If $A=\F[[x_1, \dots, x_n]]$, then  the Schouten bracket is the Buttin bracket (see \cite{CK}). 
\end{example}

\begin{lemma}(cf.\ \cite{Ki}, \cite{Ki2})
Let $A$ be a commutative associative superalgebra, $a$ an even derivation of $A$,
and $c=\sum_ib_i\wedge d_i$ an even element in $\Lambda^2(Der A)$.
Then the bracket
\begin{equation}
\{f,g\}=fa(g)-a(f)g+\sum_i(-1)^{p(f)p(b_i)}(b_i(f)d_i(g)-(-1)^{p(b_i)}
d_i(f)b_i(g))
\label{GPB}
\end{equation}
satisfies the Jacobi identity if and only if 
\begin{equation}
[a,c]_S=0 ~~\mbox{and}~~ [c,c]_S=-2a\wedge c,
\label{Sconditions}
\end{equation}
where $[\cdot,\cdot]_S$ is the Schouten bracket on $\Lambda(Der A)$.
\end{lemma}
\begin{examples}\em
$(H)$ Let $A={\cal O}(2k,n)$, $a=0$ and $c=\sum_{i,j}f_{ij}\frac{\partial}
{\partial z_i}\wedge\frac{\partial}
{\partial z_j}$, for some superskew-symmetric matrix $(f_{ij})$ over $\F$.
Then conditions (\ref{Sconditions}) hold, hence the
corresponding bracket (\ref{GPB}) defines a Poisson superalgebra
structure on $A$.

\medskip

\noindent
$(K)$ Let $A={\cal O}(2k+1,n)$, with even indeterminates $t$, $z_1, \dots,
z_{2k}$, and odd indeterminates $z_{2k+1}, \dots, z_{2k+n}$. If
$a=2\frac{\partial}{\partial t}$ and
$c=\sum_{i,j}f_{ij}\frac{\partial}
{\partial z_i}\wedge\frac{\partial}
{\partial z_j}-E\wedge\frac{\partial}{\partial t}$,
where $(f_{ij})$ is a superskew-symmetric matrix over $\F$
and $E=\sum_{i=1}^{2k+n}z_i\frac{\partial}{\partial z_i}$ is
the Euler operator,
then  conditions (\ref{Sconditions}) hold. It follows
that $A$ with bracket (\ref{GPB}) is a generalized Poisson superalgebra.

Examples \ref{P(2k,n)} and \ref{P(2k+1,n)} are particular
cases of Examples (H) and (K), respectively, when the matrix $(f_{ij})$
is invertible.
\end{examples}

\begin{example}\em
Let $A$ be a unital generalized Poisson superalgebra with bracket $\{\cdot,\cdot\}$
and derivation $D$. Then for any even invertible element 
$\varphi\in A$, formula
$$\{f,g\}^{\varphi}=\varphi^{-1}\{\varphi f, \varphi g\}$$
defines another generalized Poisson bracket on the superalgebra $A$, the
corresponding derivation being $D^\varphi(a)=D(\varphi)a+\{\varphi, a\}$.
We denote this generalized Poisson superalgebra by $A^\varphi$. The generalized
Poisson superalgebras $A$ and $A^\varphi$ are called gauge equivalent.
\end{example}

\section{Examples of Jordan superalgebras}
\begin{definition} A Jordan superalgebra is a commutative superalgebra $J$
whose product $(a,b)\mapsto a\circ b$ satisfies
the following identity:
\begin{equation}
(-1)^{p(a)p(c)}[L_{a\circ b}, L_c]+(-1)^{p(b)p(a)}[L_{b\circ c}, L_a]+
(-1)^{p(c)p(b)}[L_{c\circ a}, L_b]=0,
\label{second}
\end{equation}
where $L_a\in End(J)$
denotes the operator of left multiplication by $a$: $L_a(x)=a\circ x$.
\end{definition}

\begin{remark}\em The following important relation is proved for
Jordan superalgebras in the same way as for Jordan algebras \cite{A}:
\begin{equation}
[[L_a, L_b], L_c]=(-1)^{p(b)p(c)}L_{a\circ(c\circ b)-(a\circ c)\circ b}.
\label{relation}
\end{equation}
\end{remark}

\medskip

Now we list all the examples of simple finite-dimensional Jordan superalgebras
(see \cite[\S 1.3]{K1}). 
We denote by $gl(m,n)_+$ the full linear Jordan superalgebra, i.e.,
the Jordan superalgebra of endomorphisms of a vector superspace
$V$ of dimension
$(m|n)$, with product $a\circ b=\frac{1}{2}(ab+(-1)^{p(a)p(b)}ba)$.
 We denote by $osp(m,n)_+$ the subalgebra of $gl(m,n)_+$ consisting of
selfadjoint endomorphisms  with respect to a non-degenerate supersymmetric 
bilinear form $(\cdot,\cdot)$.
We denote by
$(m,n)_+$ the Jordan superalgebra defined on $\mathbb{F}e\oplus V$,
where $V$ is an
$(m|n)$-dimensional vector superspace with a supersymmetric bilinear form
$(\cdot,\cdot)$,
$e$ is the unit element, and the product on $V$ is defined by
$a\circ b=(a,b)e$.
We denote by $p(n)_+$ the subalgebra of $gl(n,n)_+$ consisting of
selfadjoint operators with respect to a non-degenerate odd supersymmetric
form; finally, we denote by $q(n)_+$ the subalgebra of all operators in 
$gl(n,n)_+$, commuting
with the operator $\left(\begin{array}{cc}
0 & I_n\\
I_n & 0
\end{array}\right)$.
The series $JP(0,n)$ was defined in the introduction. Finally,
for the definition
of $D_t$, $F$ and $K$ we refer to Appendix B.
The dimensions of these Jordan superalgebras $J$ are listed in 
the following table:

\bigskip

\begin{center}
\begin{tabular}{|c|c||c|c|}
\hline
$J$ & $(\dim J_{\bar{0}}|\dim J_{\bar{1}})$ & $J$ & $(\dim J_{\bar{0}}|\dim J_{\bar{1}})$\\
\hline
\vspace{0.0001cm} & \vspace{0.0001cm} & \vspace{0.0001cm} & \vspace{0.0001cm}\\
$gl(m,n)_+$ & $(m^2+n^2|2mn)$ & $D_t$ & $(2|2)$ \\
$osp(m,n)_+$ & $(\frac{1}{2}m(m+1)+\frac{1}{2}n(n-1)|mn)$   & $E$ & $(27|0)$ \\
$(m,n)_+$  & $(1+m|n)$  & $F$ & $(6|4)$ \\
$p(n)_+$ & $(n^2|n^2)$ & $K$ & $(1|2)$ \\
$q(n)_+$ & $(n^2|n^2)$ & $JP(0,n)$  & $(2^n|2^n)$\\
\hline 
\end{tabular}
\end{center}

\medskip

Next, we construct some examples of infinite-dimensional
Jordan superalgebras.
\begin{example}\label{Kdouble}\em
~Let $U$ be a commutative associative superalgebra together
with an anticommutative bilinear operation $\{\cdot,\cdot\}$,
 and consider a direct
sum of vector superspaces $J=J(U, \{\cdot,\cdot\})=U+\eta U$ where the symbol 
$\eta$ has
odd parity, i.e., the $\Z/2\Z$-grading
on $U$ is extended to $J$ by setting $J_{\bar{0}}=
U_{\bar{0}}+\eta U_{\bar{1}}$, $J_{\bar{1}}=
U_{\bar{1}}+\eta U_{\bar{0}}$. We define a multiplication on $J$ as follows:
for arbitrary elements $a, b\in U$, their product in $J$ is the
product $ab$ in $U$, and $a(\eta b)=(-1)^{p(a)}\eta(ab)$, 
$(\eta a)b=\eta(ab)$,
$(\eta a)(\eta b)=(-1)^{p(a)}\{a,b\}$. Then
$J(U, \{\cdot,\cdot\})$ is a superalgebra called the Kantor-
King-McCrimmon (KKM) double of $U$.

Let $A$ be a generalized Poisson superalgebra with bracket $\{\cdot,\cdot\}$
and derivation $D$,
and define on $A$ the following new bracket: 
$$\{f,g\}_D:=\{f,g\}-\frac{1}{2}(fD(g)-D(f)g), ~~f, g\in A.$$
Set $D':=\frac{1}{2}D$. Then the following
identities hold:
\begin{equation}
\{f,gh\}_D=\{f,g\}_Dh+(-1)^{p(f)p(g)}g\{f,h\}_D+D'(f)gh;
\label{KMC4-1}
\end{equation}
\begin{equation}
\{f,\{g,h\}_D\}_D=\{\{f,g\}_D,h\}_D+(-1)^{p(f)p(g)}
\{g,\{f,h\}_D\}_D-D'(f)\{g,h\}_D+
\label{KMC4-2}
\end{equation}

\medskip

$-(-1)^{p(f)(p(g)+p(h))}D'(g)\{h,f\}_D-
(-1)^{p(h)(p(f)+p(g))}D'(h)\{f,g\}_D.$

\medskip

\noindent
By \cite[Theorem 2]{KMCC},
the KKM double $J(A, \{\cdot,\cdot\}_D)$  is a Jordan
superalgebra. We will denote it simply by $J(A)$. This
result for Poisson superalgebras was obtained by Kantor
\cite[Theorem 3.1]{Ka2}:
notice that if $A$ is a Poisson superalgebra, then $\{\cdot,\cdot\}_D
=\{\cdot,\cdot\}$.

We set $JP(m,n):=J(P(m,n))$. If $m=0$, this is the
Kantor series of Jordan superalgebras \cite{Ka2}. All the other
Jordan superalgebras $JP(m,n)$ are inifinite-dimensional linearly
compact.
\end{example}

\begin{remark}\label{simpleunitalgen}\em
If $(A,\{\cdot,\cdot\})$ is  unital, then it is
simple if and only if $(A,\{\cdot,\cdot\}_D)$ is simple. This follows
from Remark \ref{unitalgP}.
\end{remark}

\begin{example}\label{JCK}\em Let $A$ be an associative commutative
algebra over $\F$ with a derivation $D$, and let 
$\H=\F1\oplus \H_{im}$ be the
algebra of degenerate quaternions over $\F$, i.e., 
$1$ is the unit element and
$\H_{im}=\F i\oplus\F j\oplus\F k$
is a 3-dimensional associative algebra with the product 
 $i\times j=k$, $j\times k=i$, $k\times i=j$, $i^2=j^2=k^2=0$.
Let $(\cdot, \cdot)$ be a symmetric bilinear
form on $\H_{im}$, with respect to
which $\{i, j, k\}$ is an orthonormal basis, and
extend both the product $\times$ and the form $(\cdot, \cdot)$ 
to $A\otimes \H_{im}$ by bilinearity.

The Jordan superalgebra $JCK(A,D)$ is defined as the commutative superalgebra
with even part $A\otimes\H$ and odd part 
$\eta(A\otimes\H)$,
where $\eta$ is an odd indeterminate, and the
following product (cf.\ \cite{KMZ}, \cite{MZ}):
$$\begin{array}{ll}
f\circ g=fg & \mbox{if}~ f\in A\otimes 1, ~g\in A\otimes\H;\\
f\circ g=(f,g)1 & \mbox{if}~ f,g\in A\otimes\H_{im};\\
\eta f\circ g=\eta(fg) & \mbox{if}~  g\in A\otimes 1;\\
\eta f\circ g=\eta(fD(g)) & \mbox{if}~  f\in A\otimes 1, ~g\in A\otimes\H_{im};\\
\eta f\circ g=\sqrt{-1}\,\eta(f\times g) & \mbox{if}~  f, g\in A\otimes\H_{im};\\
\eta f\circ \eta g=D(f)g-fD(g) & \mbox{if}~ f,g\in A\otimes 1;\\
\eta f\circ \eta g=fg & \mbox{if}~ f\in A\otimes\H_{im} ,~g\in A\otimes 1;\\
\eta f\circ \eta g=0 & \mbox{if}~ f,g\in A\otimes\H_{im}.
\end{array}$$
This Jordan superalgebra is simple if and only if $A$ contains no
non-trivial $D$-invariant ideals. 

If $A=\C[[x]]$ and $D=\frac{d}{dx}$, we will
denote $JCK(A,D)$ simply by $JCK$.
\end{example}

\begin{example}\label{new}\em Let $A$ be a commutative associative
superalgebra with an odd derivation $D$. Set $J=A$ with reversed
parity and define on $J$ the following product: 
$$a\circ b=aD(b)+(-1)^{p(a)p(b)}bD(a), ~~a, b\in J.$$
Then $(J, \circ)$ is a Jordan superalgebra. This is checked by a rather
long, but straightforward, calculation.

If $A={\cal O}(1,1)$ and $D=\xi\frac{\partial}{\partial x}+\frac{\partial}
{\partial\xi}$, we will denote $(J,\circ)$ by $JS$.
\end{example}

\medskip

In Section \ref{mainthms} we will prove the following theorem:

\begin{theorem}\label{main} All  linearly compact simple Jordan superalgebras
over $\F$ are, up to isomorphism, the following:
\begin{enumerate}
\item the finite-dimensional Jordan superalgebras $gl(m,n)_+$, $osp(m,2r)_+$
with $(m,r)\neq (0,1)$, $(m,2r)_+$ with $(m,r)\neq (1,0)$, $p(n)_+$ with 
$n>1$, $q(n)_+$ with $n>1$, $JP(0,n)$ with $n\geq 1$, $D_t$ with $t \in \F, t\neq 0$, 
$E$, $F$, $K$;
\item the infinite-dimensional Jordan superalgebras $JP(m,n)$ with $m\geq 1$,
\break $JCK$, and  $JS$.
\end{enumerate}
\end{theorem}

\begin{remark}\em It follows from \cite[\S 4.2.2]{K2} and the
multiplication table of the Jordan superalgebras $(m,n)_+$ and $D_t$,
that the following isomorphisms hold:
$$JP(0,1)\simeq gl(1,1)_+ ~~,~~ (1,2)_+\simeq D_1 
~~,~~ D_t\simeq D_{t^{-1}}.$$
These are the only isomorphisms of Jordan superalgebras, listed in 
Theorem \ref{main}.
\end{remark}

\begin{remark}\label{growthandsize}\em
The Jordan superalgebra $JP(m,n)$ has growth $m$ and size $2^{n+1}$;
the Jordan superalgebras $JCK$ and $JS$ have growth 1 and
sizes 8 and 2, respectively.
(For the definition of growth and size, see 
\cite[\S 1]{CantaK}.) Note that in all cases the sizes
of the even and the odd parts of $J$ are equal to $\frac{1}{2}$size$(J)$
(of course, their growths are both equal to the growth of $J$).
\end{remark}
%

\section{TKK construction}
Let $J$ be a Jordan superalgebra and  let
$P$ be the commutative bilinear  map on $J$ defined by: $P(x,y)=x\circ y$.
Notice that the $\Z/2\Z$-grading
on $J$ induces a natural $\Z/2\Z$-grading on the spaces of linear
and bilinear functions on $J$ with values in $J$.
We associate to $J$ a $\Z$-graded Lie
superalgebra $Lie(J)=\g_{-1}\oplus \g_0\oplus \g_1$ as follows
(cf.\ \cite{K1}, \cite{Ka1}, \cite{Ko}, \cite{T}). We set:
$\g_{-1}=J$,
$\g_0=\langle L_a, [L_a, L_b]~|~ a,b\in J\rangle\subset End(\g_{-1})$,
$\g_1=\langle P, [L_a, P] ~|~ a\in J\rangle$ (the bracket of 
a linear and a bilinear map is defined below),
and we define the following anticommutative brackets:
\begin{itemize}
\item[-] $[x,y]=0 ~~\mbox{for}~ x,y\in \g_{-1} ~\mbox{or}~ x,y\in\g_1;$
\item[-] $[a,x]=a(x) ~~\mbox{for}~ a\in\g_0, x\in \g_{-1};$
\item[-] $[A,x](y)=A(x,y) ~~\mbox{for}~ A\in\g_1, x,y\in \g_{-1};$
\item[-] $[a,B]=a\Box B-(-1)^{p(a)p(B)}B\Box a ~~\mbox{for}~ a\in\g_0, B\in \g_{1}$,
where \begin{itemize}
\item[i)] $(a\Box B)(x,y)=a(B(x,y))$,
\item[ii)] $(B\Box a)(x,y)=B(a(x),y)+(-1)^{p(x)p(y)}B(a(y),x)$.
\end{itemize}
\end{itemize}
$Lie(J)$ is called the Tits-Kantor-Koecher (TKK) construction for $J$.
\begin{proposition}
$Lie(J)$ is a Lie superalgebra.
\end{proposition}
{\bf Proof.} Using the definition of the bracket on $Lie(J)$ one shows that the
following relations hold:
\begin{equation}
[P,x]=L_x, ~~ x\in J;
\label{A}
\end{equation}
\begin{equation}
[[L_a,P], x]=[L_a,L_x]-L_{a\circ x}, ~~a,x\in J.
\label{*}
\end{equation}
It follows, using (\ref{*}) and (\ref{second}), that: 
\begin{equation}
[L_a,[L_b,P]]=-[L_{a\circ b},P], ~~a,b\in J;
\label{B}
\end{equation}
\begin{equation}
[[L_a,L_b],P]=0, ~~a,b\in J.
\label{C}
\end{equation}
Relations (\ref{A}) and (\ref{*}) show that $[\g_{-1}, \g_1]\subset\g_0$.
Besides, a direct computation shows that:
\begin{equation}
[a,[b,C]]=[[a,b],C]+(-1)^{p(a)p(b)}[b,[a,C]], ~~\mbox{for}~
a,b\in\g_0, C\in\g_1.
\label{added}
\end{equation}
It follows from (\ref{added}), (\ref{C}) and (\ref{relation}), that 
\begin{equation}
[[L_a,L_b],[L_c,P]]=(-1)^{p(c)p(b)}[L_{a\circ(c\circ b)-(a\circ c)\circ b},P].
\label{D}
\end{equation}
Relations (\ref{B}), (\ref{C}) and (\ref{D}) show that 
$[\g_0,\g_1]\subset \g_1$.
Finally, using identity (\ref{relation}),  one verifies 
that $\g_0$ is a subalgebra of $Lie(J)$. Therefore $Lie(J)$ is closed
under the defined bracket.
Besides, a direct computation, together with relations
(\ref{A})-(\ref{C}), shows that 
the Jacobi identity is satisfied.  
\hfill$\Box$ 

\begin{definition}\label{shortgrading}
A short grading of a Lie superalgebra $\g$ is a $\Z$-grading
of the form $\g=\g_{-1}\oplus \g_0\oplus \g_1$.
\end{definition}

\begin{definition}\label{short}
A short subalgebra of a Lie superalgebra $\g$
is an $sl_2$ subalgebra spanned by even elements $e$, $h$, $f$ such that:
\begin{itemize}
\item[$(i)$] the eigenspace decomposition of $ad \, h$ defines a short
grading of $\g$;
\item[$(ii)$] $[h,e]=-e$, $[h,f]=f$, $[e,f]=h$.
\end{itemize}
\end{definition}

\begin{example}\label{triple}\em Let $J$ be a Jordan
superalgebra with unit element $e$ and let $Lie(J)$ be the Lie
superalgebra  associated to $J$. Consider the elements
$h_J:=-L_e$ and
$f_J=P$ in $Lie(J)$. Then the subalgebra $\mathfrak{a}_J$ spanned by 
$e$, $f_J$ and
$h_J$ is a short subalgebra of $Lie(J)$.
\end{example}

\begin{lemma}\label{exp} The adjoint representation of a short subalgebra
  $\mathfrak{a}$ of a Lie superalgebra $\g$
  exponentiates to the action of $SL_2$ by continuous automorphisms on $\g$.
\end{lemma}
{{\bf Proof.} By Definition \ref{short}, $\g=\g_{-1}\oplus\g_0\oplus\g_1$
where $\g_i=\{x\in\g~|~[h,x]=ix\}$. It follows that $(ad\, e)^3=0$ and
$(ad\, f)^3=0$. 
  Elements $e$ and $f$  thus exponentiate to continuous automorphisms
  of $\g$ since $\exp(ad\, e)=1+ad\, e+\frac{1}{2}(ad\, e)^2$ and
  $\exp(ad\, f)=1+ad\, f+\frac{1}{2}(ad\, f)^2$. Since $\a$ is generated by 
$e$ and $f$, the thesis follows.
\hfill$\Box$

\begin{definition}\label{minimal} 
A $\Z$-graded Lie superalgebra  $\g=\g_{-1}\oplus \g_0\oplus \g_1$ 
(and the corresponding grading) is called minimal if
any non-trivial ideal $I$ of $\g$ intersects $\g_{-1}$ non-trivially,
i.e., $I\cap\g_{-1}$ is neither 0 nor $\g_{-1}$.

If $\g$ is a Lie superalgebra and
$\a$ is a short subalgebra of $\g$ such that
the corresponding short grading $\g=\g_{-1}\oplus \g_0\oplus \g_1$ is
minimal, then 
the pair $(\g, \a)$ is called a minimal pair.
\end{definition}

\begin{lemma}\label{equivcond}
A $\Z$-graded Lie superalgebra
$\g=\g_{-1}\oplus\g_0\oplus\g_1$ is minimal if and only if 
the following three conditions hold:
\begin{equation}
\mbox{if}~~ [a,\g_{-1}]=0 ~\mbox{for some}~ a\in \g_0\oplus\g_1, 
~\mbox{then}~ a=0~~(\mbox{transitivity});
\label{(i)}
\end{equation}
\begin{equation}
[\g_{-1},\g_1]=\g_0;
\label{(ii)}
\end{equation}
\begin{equation}
[\g_0, \g_1]=\g_1.
\label{(iii)}
\end{equation}
\end{lemma}
{\bf Proof.} Let $\g=\g_{-1}\oplus\g_0\oplus\g_1$ be a minimal
grading and let $A=\{a\in\g_0\oplus\g_1~|~[a,\g_{-1}]=0\}$. If $A\neq 0$,
then $I=\sum_{n\geq 0}(ad\, \g_0)^n(A)+\sum_{n\geq 0}(ad\, \g_0)^n([\g_1,A])
\subset \g_0\oplus\g_1$ is a non-trivial ideal of $\g$ with
$I\cap\g_{-1}=0$, hence contradicting the
minimality of the grading. Transitivity thus follows. Likewise,
if $[\g_{-1},\g_1]\subsetneq\g_0$, then $I=\g_{-1}\oplus [\g_{-1},\g_1]
\oplus \g_1$ is a non-trivial ideal of $\g$ 
with $I\cap\g_{-1}=\g_{-1}$, thus contradicting the hypotheses.
Finally, if $[\g_0,\g_1]\subsetneq\g_1$, then $I=\g_{-1}\oplus\g_0\oplus
[\g_0,\g_1]$ is a non-trivial ideal of $\g$ such that
$I\cap\g_{-1}=\g_{-1}$, thus contradicting the minimality
of the grading.

Conversely, assume that conditions (\ref{(i)}), (\ref{(ii)}) and (\ref{(iii)})
 hold
and let $I$ be a non-trivial ideal of $\g$. Then 
$I\cap\g_{-1}\neq 0$ by condition (\ref{(i)}), and $I\cap\g_{-1}\neq \g_{-1}$
by conditions (\ref{(ii)}) and (\ref{(iii)}).
\hfill$\Box$

\begin{lemma}\label{properties}
Let $J$ be a Jordan superalgebra and let $Lie(J)=
\g_{-1}\oplus\g_0\oplus\g_1$ be the corresponding Lie superalgebra. Then
conditions (\ref{(i)}) and (\ref{(ii)}) hold.
If, in addition, $J$ is unital, then also condition (\ref{(iii)}) holds,
hence $(Lie(J), \a_J)$ is a minimal pair.
\end{lemma}
{\bf Proof.} 
(\ref{(i)}) follows immediately from the definition of the bracket in
$Lie(J)$ and (\ref{(ii)}) follows from properties (\ref{A}) and (\ref{*}).
Finally, if $e$ is the unit element of $J$, then $[L_e,P]=-P$ which
proves (\ref{(iii)}). Hence, by Lemma \ref{equivcond},
if $J$ is unital, then $(Lie(J), \a_J)$ is a minimal pair. \hfill$\Box$

\begin{remark}\em If a short grading $\g=\g_{-1}\oplus\g_0\oplus\g_1$
comes from a short subalgebra,
then it is minimal if and only if $\g_0$ contains no non-zero ideals of $\g$. 
\end{remark}

\begin{proposition}\label{simplicity}
A unital Jordan superalgebra $J$ is  simple if and only if $Lie(J)$ is a simple
Lie superalgebra.
\end{proposition}
{\bf Proof.} By Lemma \ref{properties}, $Lie(J)=\g_{-1}\oplus \g_0\oplus\g_1$
satisfies (\ref{(i)}), (\ref{(ii)}), (\ref{(iii)}). Suppose that $J$ is simple. Then
$\g_{-1}$ is an irreducible $\g_0$-module.
Let $I\neq 0$ be an ideal of $Lie(J)$. Then, by 
(\ref{(i)}), $I\cap\g_{-1}\neq 0$,
hence, by the irreducibility of $\g_{-1}$, $I\supset\g_{-1}$. 
Then the simplicity of $Lie(J)$ follows from
(\ref{(ii)}) and (\ref{(iii)}).

Conversely, if $Lie(J)=\g_{-1}\oplus\g_0\oplus\g_1$ is simple, 
then $\g_{-1}$ is an irreducible $\g_0$-module. Indeed, if
$\h_{-1}$ is a proper $\g_0$-submodule of $\g_{-1}$,
then $\h_{-1}\oplus [\h_{-1}, \g_1]\oplus[[\h_{-1}, \g_1], \g_1]$
is a proper ideal of $Lie(J)$. If $Y\neq 0$ is
an ideal of $J=\g_{-1}$, then  $Y$ is
 $\g_0$-stable, hence $Y=J$ by the irreducibility of $\g_{-1}$.
 \hfill$\Box$ 

\bigskip

If $J$ is a non-unital Jordan superalgebra, we define $\tilde{J}=\F 1\oplus J$
and we extend the product of $J$ to a product on $\tilde{J}$ by
setting $1\circ x=x$. Then the following proposition holds:

\begin{proposition}\label{nonunital} If $J$ is a simple non-unital Jordan superalgebra, 
then $Lie(\tilde{J})$ $\simeq S\rtimes sl_2$, where $S$ is a
simple Lie superalgebra and $sl_2$ acts
on $S$ by outer derivations. If $J$ is a non-unital Jordan superalgebra such 
that
$Lie(\tilde{J})$ is as above, then $J$ is
simple. 
\end{proposition}
{\bf Proof.} By relations (\ref{A}) and (\ref{B}),
elements $1$, $-L_1$ and $P$ span a short subalgebra $\a$ of
$Lie(\tilde{J})$. Let $Lie(\tilde{J})=Lie(\tilde{J})_{-1}\oplus
Lie(\tilde{J})_0\oplus Lie(\tilde{J})_1$ be the corresponding
short grading. By relations (\ref{*})-(\ref{C}), $Lie(\tilde{J})=S\rtimes \a$,
where  $S=S_{-1}\oplus S_0\oplus S_1$, $S_{-1}=J$,
$S_1=\langle [P,L_a] ~|~ a\in J\rangle$, and $S_0$ is the subalgebra
of $Lie(\tilde{J})_0$ spanned by $L_a$ and $[L_a, L_b]$ with $a,b$ in $J$.
Here $\a$ acts on $S$ by outer derivations.
The element 
$\left(\begin{array}{cc}
0 & -1\\
1 & 0
\end{array}
\right)\in SL_2$
exchanges $Lie(\tilde{J})_{-1}$  with $Lie(\tilde{J})_{1}$, 
hence it exchanges $S_{-1}$ with $S_1$, 
but the $S_0$-module $S_{-1}$ is irreducible
since $J$ is simple, hence $[S_0, S_{-1}]=S_{-1}$
(by definition, a simple algebra is an algebra with 
non-trivial product and no non-trivial ideals). It follows that
$[S_0, S_1]=S_{1}$. Now the simplicity of $S$
is proved as the simplicity of $Lie(J)$ 
in Proposition \ref{simplicity}.
The second part is proved as 
in Proposition \ref{simplicity} too.
\hfill$\Box$

\begin{corollary}\label{LieJmin}
Let $J$ be a Jordan superalgebra. Then $Lie(J)=\g_{-1}\oplus\g_0\oplus\g_1$
is a minimal $\Z$-graded superalgebra.
\end{corollary}
{\bf Proof.} If $J$ is unital then the statement holds by
Lemma \ref{properties}. Now suppose that $J$ is not unital. 
By Lemma \ref{equivcond},
$Lie(J)=\g_{-1}\oplus\g_0\oplus\g_1$
is a minimal $\Z$-graded superalgebra if and only if it satisfies
properties (\ref{(i)}), (\ref{(ii)}), (\ref{(iii)}). Properties 
(\ref{(i)}), (\ref{(ii)}) hold by Lemma \ref{properties}, and
property (\ref{(iii)}) holds by the proof of Proposition \ref{nonunital}.
\hfill$\Box$

\bigskip

Let $\cal{J}$ denote the category of unital Jordan superalgebras in which
morphisms are usual Jordan superalgebra homomorphisms, and let
${\cal L}$ denote the category of minimal
 pairs $(\g,\a)$, where $\g$ is a Lie
superalgebra, $\a$ is a short subalgebra of
$\g$, and a morphism $\varphi$ from a pair $(\g,\a)$ to a pair
$(\g', \mathfrak{b})$ is a Lie superalgebra homomorphism $\varphi: \g
\rightarrow \g'$ such that $\varphi(\mathfrak{a})=\mathfrak{b}$.
According to Lemma \ref{properties}, we construct a functor 
$\cal{F}$ from ${\cal J}$ to ${\cal L}$ by associating to
every unital Jordan superalgebra $J$ the pair 
$(Lie(J),\a_J)$,
and to every homomorphism $\varphi: J\rightarrow J'$ of unital Jordan
superalgebras 
 the map $\varphi_{\cal F}:
Lie(J) \rightarrow Lie(J')$ defined as follows:
$$x\mapsto \varphi(x), ~ L_x\mapsto L_{\varphi(x)} ~~\mbox{for every}~ x\in J;
~~P\mapsto P'$$
where $P'$ is the bilinear  map on $J'$ defined by: $P'(a,b)=a\circ b$ for
$a,b\in J'$.
One checks that $\varphi_{\cal F}$ is a Lie algebra homomorphism
mapping $\a_J$ to $\a_{J'}$.
We will show that ${\cal F}$ defines  an equivalence of 
categories. In order to construct the inverse functor
to ${\cal F}$
we need the following result \cite{T}, \cite[Lemma 4]{K1}:
\begin{lemma}\label{lemma4}
Let $\g=\g_{-1}\oplus\g_0\oplus\g_1$ be a $\Z$-graded Lie superalgebra
and let $p\in(\g_1)_{\bar{0}}$. Set
\begin{equation}
x\circ y:=[[p,x],y].
\label{Jproduct}
\end{equation}
Then the $\Z/2\Z$-graded space $\g_{-1}$ with product
(\ref{Jproduct}) is a Jordan superalgebra.
\end{lemma}

\begin{remark}\label{inversefunctor}\em
Let $\g$ be a Lie superalgebra containing a short subalgebra $\mathfrak{a}$
with standard generators $\{e,h,f\}$ as in Definition \ref{short}.
Consider the corresponding short grading:
 $\g=\g_{-1}\oplus\g_0\oplus\g_1$ where $\g_i=\{x\in\g~|~[h,x]=ix\}$.
Define on $J(\g,\a):=\g_{-1}$ product (\ref{Jproduct}), where $p=f$.
Then, by Lemma \ref{lemma4},
$J(\g,\a)$ is a Jordan superalgebra and $e$ is its unit element.
\end{remark}

\begin{theorem}\label{equivalence}
The functor ${\cal F}$ is an equivalence of categories.  
\end{theorem}
{\bf Proof.} In order to define the inverse functor
to ${\cal F}$ we  associate to the minimal pair $(\g, \a)$,
where $\g$ is a Lie
superalgebra and $\a$ is a 
short subalgebra of $\g$, the Jordan superalgebra $J(\g,\a)$ defined
in Remark \ref{inversefunctor}, and to every Lie superalgebra
homomorphism $\psi: (\g,\a) \rightarrow (\g',\mathfrak{b})$, the Jordan
superalgebra homomorphism $\psi_{\Phi}:J(\g,\a)\rightarrow J(\g',\mathfrak{b})$, 
$\psi_{\Phi}=\psi_{|\g_{-1}}$.
Notice that, since $\psi$
is a Lie superalgebra homomorphism mapping $\a$ to $\mathfrak{b}$,
due to Lemma \ref{exp}, we may assume that the corresponding triples 
$\{e,h,f\}$ are mapped to each 
other,
hence, by (\ref{Jproduct}), $\psi$
induces a homomorphism of Jordan superalgebras. Therefore 
we have defined a functor from ${\cal L}$ to ${\cal J}$. Let us denote
it by $\Phi$.

Let $(\g,\a)$ be a minimal pair, with $\a=\langle e,f,h\rangle$ and
$\g=\g_{-1}\oplus\g_0\oplus\g_1$ the corresponding minimal grading,
and let $(Lie(J), \a_J)$ be the pair associated to $J(\g,\a)$ by
${\cal F}$. 
From the representation theory of $\a\simeq sl_2$, applied to its adjoint 
action  on $\g$, it follows that:
$[e,\g_0]=\g_{-1}$, $[f,\g_0]=\g_1$ and $\g_0 = [\g_{-1},f]+ \g_0^f$,
where $\g_0^f$ is the centralizer of $f$ in $\g_0$.
Since, by the minimality of $(\g,\a)$, $\g_0=[\g_1, \g_{-1}]$,
it follows that $\g_0
=[[\g_0,f],\g_{-1}]=[[[\g_{-1},f],f],\g_{-1}]=[[\g_{-1},f],[\g_{-1},f]] +
[\g_{-1},f]$, hence $(Lie(J), \a_J)=(\g, \a)$. Besides,
if $\psi: (\g, \a) \rightarrow (\g', \b)$ is a Lie superalgebra
homomorphism such that $\psi(\a)=\b$, then $(\psi_\Phi)_{\cal F}=\psi$.

Conversely, let $J$ be a Jordan superalgebra with unit element $e$. 
Consider the
pair $(Lie(J), \a_J)$ associated to $J$ by ${\cal F}$ and
the unital Jordan superalgebra $J(Lie(J),$ $\a_J)$ associated to
$(Lie(J), \a_J)$ by $\Phi$. By construction, if $a$ and $b$ lie
in $J(Lie(J),$ $\a_J)$, their product is $a\circ b=[[[P, L_e],a],b]$.
By (\ref{*}), this is the same as the product of $a$ and $b$ in $J$. 
Besides, if $\varphi: J \rightarrow J'$
is a homomorphism of unital Jordan superalgebras, and $\varphi_{\cal F}:
(Lie(J), \a_J) \rightarrow (Lie(J'), \a_{J'})$
is the corresponding homomorphism of pairs in ${\cal L}$,
with $Lie(J)=\g_{-1}\oplus\g_0\oplus\g_1$, then
${\varphi_{\cal F}}_{|\g_{-1}}=\varphi$.
We have thus shown that $\Phi$ is the inverse
functor to ${\cal F}$.
\hfill$\Box$

%

\bigskip

Let ${\cal J}_{lc}$ denote the subcategory of ${\cal J}$ consisting of
 unital linearly compact
 Jordan superalgebras in which
morphisms are continuous Jordan superalgebra homomorphisms, and let
${\cal L}_{lc}$ denote the subcategory of ${\cal L}$ consisting of
minimal pairs $(\g,\a)$, where $\g$ is a linearly compact Lie
superalgebra, and morphisms from a pair $(\g,\a)$ to a pair
$(\g', \mathfrak{b})$ are continuous Lie superalgebra homomorphisms 
$\varphi: \g
\rightarrow \g'$ such that $\varphi(\mathfrak{a})=\mathfrak{b}$.

The TKK construction carries over to the case of a linearly
compact $J$ almost verbatim, except that in the construction of
$\g_0$ in $Lie(J)$ one takes the closure of the corresponding 
span.
\begin{proposition}\label{lc}
The functor ${\cal F}$ is an equivalence of the categories
${\cal J}_{lc}$ and ${\cal L}_{lc}$.
\end{proposition}
{\bf Proof.} ${\cal F}$ is a functor from  
${\cal J}_{lc}$ to ${\cal L}_{lc}$ due to Lemma \ref{lincomp}.
Besides, functor $\Phi$ maps ${\cal L}_{lc}$ to
${\cal J}_{lc}$ since if 
$\g$ is a linearly compact Lie superalgebra and $\a$ is a short subalgebra
of $\g$ inducing on $\g$ the grading $\g=\g_{-1}\oplus\g_0\oplus\g_1$, then
$\g_{-1}$ is a closed
subspace of $\g$, hence it is linearly compact.
By construction, continuous homomorphisms are mapped to continuous 
homomorphisms. 
\hfill$\Box$

\section{Classification of short gradings of simple linearly compact Lie
superalgebras}
Let $L$ be a simple finite-dimensional Lie algebra and choose a
vertex $i$ of its Dynkin diagram. Consider the $\Z$-grading of $L$
defined by setting $\deg(e_i)=1$, $\deg(f_i)=-1$ and the degrees of
all other Chevalley generators of $L$ equal to 0. We will denote 
such a grading by the Dynkin diagram of $L$ where the $i$-th vertex
is labelled by $+$. The following well known
 theorem describes all gradings defined by short subalgebras
of all simple finite-dimensional Lie algebras:

\begin{theorem}\label{sg} A complete list, up to
isomorphism, of short gradings defined by short subalgebras of all
simple finite-dimensional Lie algebras is
as follows:

$\begin{array}{ccc}
& & \\
A_{2n-1} & (n\geq 1): & 
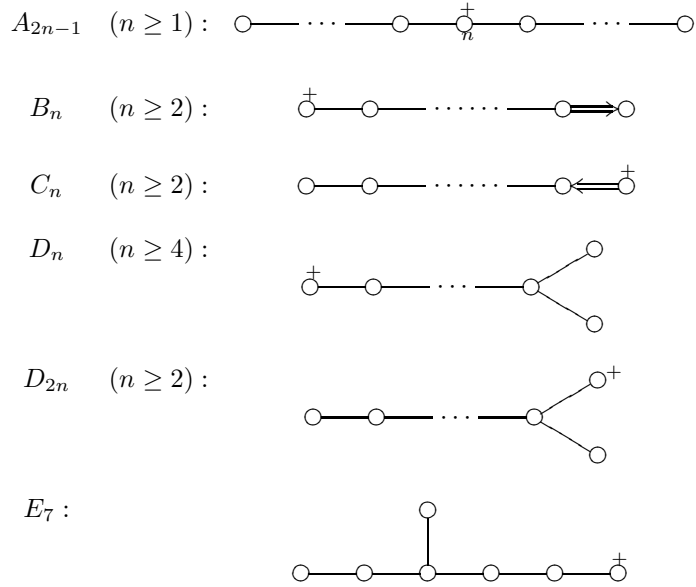
\xymatrix@R=6pt@C=18pt{
*+[o][F]{}\ar@{-}[r] &\cdots &*+[o][F]{}\ar@{-}[r]\ar@{-}[l]  & 
*+[o][F]{}\ar@{-}^{+~~~~~~}_{n~~~~~~}[r]  & *+[o][F]{}\ar@{-}[r]
&\cdots&  *+[o][F]{}\ar@{-}[l]
}\\
& & \\
B_n & (n\geq 2): & \xymatrix@R=6pt@C=18pt{
*+[o][F]{}\ar@{-}[r]^{+~~~~~~} &  *+[o][F]{}\ar@{-}[r] & \cdots\cdots  &
*+[o][F]{}\ar@{-}[l]& *+[o][F]{}\ar@{<=}[l] \\
}\\
& & \\
C_n & (n\geq 2): & ~\,\xymatrix@R=6pt@C=18pt{
*+[o][F]{}\ar@{-}[r] &  *+[o][F]{}\ar@{-}[r] & \cdots\cdots  &
*+[o][F]{}\ar@{-}[l]\ar@{<=}[r]^{~~~~~~~+}& *+[o][F]{} \\
}\\
& & \\
D_n & (n\geq 4): & \xymatrix@R=6pt@C=18pt{
&&&& *+[o][F]{}\ar@{-}[ld]\\
*+[o][F]{}\ar@{-}[r]^{+~~~~~~} &  *+[o][F]{}\ar@{-}[r] & \cdots  &
 *+[o][F]{}\ar@{-}[l] \\
&&&& *+[o][F]{}\ar@{-}[lu]}~~~\\
& & \\
D_{2n} & (n\geq 2): & \xymatrix@R=6pt@C=18pt{
&&&& *+[o][F]{}\ar@{-}[ld]\\
*+[o][F]{}\ar@{-}[r] &  *+[o][F]{}\ar@{-}[r] & \cdots  &
 *+[o][F]{}\ar@{-}[l] \\
&&&& *+[o][F]{}\ar@{-}[lu]}^{+}\\
& & \\
E_7: & &
\xymatrix@R=6pt@C=18pt{
&&  *+[o][F]{}\ar@{-}[d]\\
& & *+-[][F-]{}\ar@{-}[d] &\\
*+[o][F]{}\ar@{-}[r]&
*+[o][F]{}\ar@{-}[r]  & *+[o][F]{}\ar@{-}[r] &*+[o][F]{}\ar@{-}[r]  & *+[o][F]{}\ar@{-}[r]^{~~~~~~~+}  & *+[o][F]{}\\
}
\end{array}$
\end{theorem} 

\medskip

\noindent
{\bf Proof.} Any $\Z$-grading $L=\oplus_{j\in\Z}L_j$ of a simple
finite-dimensional Lie algebra $L$ is defined by letting $\deg(e_i)=
-\deg(f_i)=k_i\in\Z_+$, where $e_i$, $f_i$ are the Chevalley
generators of $L$. Let $\sum_ia_i\alpha_i$ be the highest root of $L$.
It is clear that a $\Z$-grading is short if and only if all $k_i$'s are
$0$, except for $k_s=1$, such that $a_s=1$. Next, if
a short grading comes from a short subalgebra $\a=\langle e, h, f\rangle$,
then: 
\begin{equation}
(L_0^e,h)=0,
\label{victorparis}
\end{equation}
where $(\cdot,\cdot)$ is the Killing form of $L$ and $L_0^e$
denotes the centralizer of $e$ in $L_0$. Indeed: $(L_0^e, [e,f])=
([L_0^e,e],f)=0$. Conversely, if (\ref{victorparis}) holds,
then $h\in[e,L_1]$, since $[e, L_1]\perp L_0^e$ and
$L_0=L_0^e+[e,L_1]$ (by the representation theory of $sl_2$).
Thus (\ref{victorparis}) is a necessary and sufficient condition for 
a short grading to come from a short subalgebra, and this 
condition (along with $a_s=1$) produces the above list.
\hfill$\Box$

\begin{remark}\em We denote by 
$\a(sl_{2n})$,
$\a'(so_{2n+1})$, $\a(sp_{2n})$, $\a'(so_{2n})$, 
$\a(so_{4n})$, and $\a(E_7)$,
respectively, the short subalgebras of the simple Lie algebras 
corresponding to the diagrams listed in Theorem \ref{sg}.
If for the definition of the classical Lie algebras  $so_{2n+1}$, 
$sp_{2n}$, and $so_{2n}$
we choose the bilinear forms with matrices
$$\left(\begin{array}{ccc}
1 & 0 & 0\\
0 & 0 & I_n\\
0 & I_n & 0
\end{array}\right),
~~~\left(\begin{array}{cc}
0 & I_n\\
-I_n & 0
\end{array}\right), 
~~~\left(\begin{array}{cc}
0 & I_n\\
I_n & 0
\end{array}\right),$$
respectively, in some basis,
then  the
elements $h$ in the corresponding short subalgebras are the following:
\begin{enumerate}
\item[] $\a(sl_{2n})$:  ~~~~~~~$h=\frac{1}{2}diag(I_n, -I_n)$;
\item[] $\a'(so_{2n+1})$: ~~\,$h=diag(0,1,0, \dots,0, -1,0, \dots,0)$;
\item[] $\a(sp_{2n})$: ~~~~~~$h=\frac{1}{2}diag(I_n, -I_n)$;
\item[] $\a'(so_{2n})$: ~~~~~$h=diag(1,0, \dots,0, -1,0, \dots,0)$;
\item[] $\a(so_{4n})$: ~~~~~\,$h=\frac{1}{2}diag(I_{2n}, -I_{2n})$.
\end{enumerate}
\end{remark}

\begin{theorem}\label{fdshortsa} A complete list, up to isomorphism,
 of pairs $(S,\a)$
where $S$ is a simple linearly compact 
Lie superalgebra, and
$\a$ is a short subalgebra of $Der S$, is as follows:
\begin{enumerate}
\item $S=sl(2m,2n)/\delta_{m,n}\F I_{2m+2n}$, $\a$ is the diagonal subalgebra of $\a(sl_{2m})+
\a(sl_{2n})$;
\item $S=osp(m,2r)$ with $m\geq 3$ and $(m,r)\neq (4,0)$, $\a=\a'(so_m)$;
\item $S=osp(4l,2r)$ with $(l,r)\neq (1,0)$, $\a$ the diagonal subalgebra of $\a(so_{4l})+
\a(sp_{2r})$;
\item $S=p(2k)$ with $k>1$, $\a=\a(sl_{2k})$;
\item $S=q(2k)$ with $k>1$, $\a=\a(sl_{2k})$;
\item $S=E_7$, $\a=\a(E_7)$;
\item $S=F(4)$, $\a$ is the diagonal subalgebra of $\a(sl_2)+\a(so_7)$;
\item $S=D(2,1;\alpha)$, $\a$ is the diagonal subalgebra of
$\a(sl_2)+\a(sl_2)$;
\item $S=sl(2,2)/\F I_4$, $\a$ is the subalgebra of outer derivations
of $sl(2,2)/\F I_4$ (see \cite[Proposition 5.1.2(e)]{K2});
\item $S=H(2k,n)$ with $n\geq 3$ and $(k,n)\neq (0,3)$, $\a=\langle \xi_n\xi_{n-1},
\xi_{n-2}\xi_n, \xi_{n-2}\xi_{n-1}\rangle$;
\item $S=K(2k+1,n)$ with $n\geq 3$, $\a=\langle \xi_n\xi_{n-1},
\xi_{n-2}\xi_n, \xi_{n-2}\xi_{n-1}\rangle$;
\item $S=E(1,6)$ and $\a=\langle \xi_5\xi_6,
\xi_4\xi_6, \xi_4\xi_5\rangle$;
\item $S=S(1,2)$, $\a$ is the subalgebra of outer derivations 
of $S(1,2)$ (see Example \ref{JS(1,1)} below).
\end{enumerate}
\end{theorem}
{\bf Proof.} 
Let $S$ be a simple linearly compact Lie superalgebra, let $Der S$ be the
Lie superalgebra of its continuous derivations (recall that it is linearly 
compact artinian semisimple \cite{CantaK}), and let $\a$ be a short subalgebra 
with the basis $e$, $h$ and $f$ as in Definition \ref{short}, of the Lie 
superalgebra $Der S$. Then it is a
short subalgebra of a maximal reductive 
subalgebra of $Der S$.
Recall that a finite-dimensional even subalgebra is called reductive if its 
adjoint representation on the whole Lie superalgebra decomposes into a product 
of finite-dimensional irreducible submodules, that all maximal reductive 
subalgebras of an artinian semisimple linearly compact Lie superalgebra 
($Der S$ is the superalgebra in question) are
conjugate to each other, and any reductive subalgebra is contained in one of 
them (see \cite[Theorem 1.7]{CantaK}; this theorem is stated for tori,
but its proof works for reductive subalgebras as well).

We will go over the list of all finite- and infinite-dimensional
simple linearly compact Lie superalgebras $S$ \cite{K2}, \cite{K3}.
Let us first assume that $\a\subset S$.
If $S=sl(r,s)$ with $r\neq s$, 
then, since $S_{\bar{0}}=sl_r\oplus sl_s\oplus\F$, 
by Theorem \ref{sg}, $r=2m$, $s=2n$, and 
$h$ is one of the following: $h_1=1/2 diag(I_m, -I_m, I_n, -I_n)$,
$h_2=1/2 diag(I_m,-I_m,0,0)$, $h_3=1/2 diag(0,0,I_n,-I_n)$. 
Using the description of
the odd roots of $S$ given in \cite[\S 2.5]{K2}, one sees that
$ad\,h_2$ and $ad\,h_3$ have non-integers eigenvalues on $S_{\bar{1}}$,
hence $h=h_1$. The same argument takes care of all other finite-dimensional
$S$ with $S_{\bar{0}}$ reductive.
For example,
if $S=p(n)$ or $S=q(n)$, then $S_{\bar{0}}=sl_{n}$, hence $S$ has, up
to conjugation, at most one short subalgebra, and it exists if and only if 
$n$ is even. 

Now consider $S=W(m,n)$ or $S=S(m,n)$ with the $\Z$-grading 
defined by:
$\deg(x_i)=-\frac{\partial}{\partial x_i}=1$; $\deg(\xi_j)=
-\deg(\frac{\partial}{\partial\xi_j})=1$.
With this grading $S=\prod_{j\geq -1}S_{j}$ where $S_0\simeq gl(m,n)$
if $S=W(m,n)$ and $S_0\simeq sl(m,n)$ if $S=S(m,n)$.
Then, 
$\a$ is conjugate to a short subalgebra of 
$sl(m,n)$. Since  $S_{-1}$ is isomorphic
to the standard $S_0$-module, from the description of short
subalgebras in $sl(m,n)$ it follows that $h$ has non-integer eigenvalues
on $S_{-1}$. Hence $W(m,n)$ and $S(m,n)$ have no short subalgebras.
Similar arguments show that $\tilde{S}(0,n)$,
$HO(n,n)$, $SHO(n,n)$, $KO(n,n+1)$, $SKO(n,n+1;\beta)$,
$SHO^\sim(n,n)$ and $SKO^\sim(n,n+1)$ 
have no short subalgebras.

Let $S=H(2k,n)$ and let $\a=\langle e,h,f\rangle$ be a short subalgebra of $S$.
 Consider $S$ with
the $\Z$-grading defined by $S_j=\{P(x, \xi)~|~\deg P=j+2\}$.
Then $S=\prod_{j\geq -1}S_j$ where $S_0\simeq osp (n,2k)$ and $S_{-1}$ is
isomorphic to the standard $osp(n,2k)$-module. Therefore
$\a$ is conjugate to a short subalgebra of $osp(n,2k)$ such that
$ad\,h$ has integer eigenvalues on $S_{-1}$. It follows that $n\geq 3$ and
$\a$ is conjugate to $\a'(so_n)$, i.e., up to conjugation,
$\a=\langle \xi_n\xi_{n-1}, \xi_{n-2}\xi_n,$ $\xi_{n-2}\xi_{n-1}\rangle$.
Similar arguments show that $K(m,n)$
with $n\geq 3$, and
$E(1,6)$ have, up to conjugation, a unique
short subalgebra, and that $K(m,n)$ has no short subalgebras if $n\leq 2$.

If $S=E(5,10)$, then the maximal reductive subalgebra of $S$ is
isomorphic to $sl_5$, hence, by Theorem \ref{sg}, $S$ has no short subalgebras.

We recall that a $\Z$-grading of a Lie superalgebra $L$ is
equivalent to the choice of an ad-diagonalizable element $\t$
with integer eigenvalues in $Der L$. If $L$ is an artinian semisimple
linearly compact Lie superalgebra, then
all maximal tori are conjugate by an inner automorphism of $L$ 
\cite[Theorem 1.7]{CantaK}. It follows that,
up to conjugation by an inner automorphism of $L$, the $\Z$-gradings
of $L$ are parameterized by the set $P^\vee=\{\t\in T
~|~ad(t) ~\mbox{has only integer eigenvalues}\}$, where $T$ is a fixed
maximal torus of $L$.

Let $S=E(3,6)$. It follows from the above discussion and 
\cite{CK} that all
$\Z$-gradings of $S$ are parameterized by quadruples $(a_1, a_2,
a_3, \varepsilon)$ such that $\varepsilon+\frac{1}{2}\sum_{i=1}^3 a_i\in\Z$,
where:
\begin{equation}
\deg(x_i)=-\deg(\frac{\partial}{\partial x_i})=a_i\in\Z, ~~\deg(d)=-\frac{1}{2}\sum_{i=1}^3 a_i,
\label{a}
\end{equation}
\begin{equation}
\deg(e_1)=-\deg(e_2)=\varepsilon\in\frac{1}{2}\Z,
\label{b}
\end{equation}
\begin{equation}
\deg(E)=-\deg(F)=2\varepsilon.
\label{c}
\end{equation}
It is immediate to see that a necessary condition for such a grading to be
short (in fact, to be finite) is that $a_i=0$ for every $i=1,2,3$. But this implies $\varepsilon\in\Z$,
and (\ref{c}) shows that the resulting grading cannot be short. It follows
that $E(3,6)$ has no short subalgebras. Similar arguments, using
the explicit description of $\Z$-gradings \cite{CantaK}, show
 that the Lie superalgebras
$E(3,8)$ and $E(4,4)$ have no short gradings, hence no
short subalgebras.

According to \cite[Proposition 5.1.2(e)]{K2}, \cite[Proposition 6.1]{K3}
and its corrected version \cite[Proposition 1.8]{CantaK},
if $S$ is a simple linearly compact Lie superalgebra such that
$Der S$ contains a copy of
$sl_2$ of outer derivations of $S$, then $S$ is isomorphic to
$sl(2,2)/\F I_4$, $S(1,2)$, $SHO(3,3)$, or $SKO(2,3;1)$.
One can show directly  that 
these $sl_2$'s in $S=sl(2,2)/\F I_4$ and  $S=S(1,2)$
(see Example \ref{JS(1,1)} below) are short subalgebras
of $Der S$. On the contrary, if  $S=SHO(3,3)$ or
$S=SKO(2,3;1)$, then $Der S$ has no short gradings, hence it has no short
subalgebras. Indeed, every $\Z$-grading of $SHO(3,3)$ is, up to conjugation,
the $\Z$-grading defined by:
$$\deg(x_i)=-\deg(\frac{\partial}{\partial x_i})=a_i\in\Z; ~~\deg(\xi_j)=
-\deg(\frac{\partial}{\partial\xi_j})=b_j\in\Z$$
for some $a_i$'s, $b_j$'s such that $a_i+b_i=c\in\F$. A necessary
condition for such a grading to be short is that $a_i=0$ for every $i=1,2,3$,
hence $b_1=b_2=b_3=1$ or $-1$. But the induced grading of
$Der S$ is not short (for a description of $S$ and $Der S$ see \cite{K3}).
Similar arguments show that $Der SKO(2,3;1)$ has no short subalgebras.
\hfill$\Box$

\begin{remark}\label{finitedimcases}\em 
Let $(S,\a)$ be one of the pairs
listed in Theorem \ref{fdshortsa}, with $S$  
a simple finite-dimensional Lie superalgebra different from $H(0,n)$ (which is
treated in Example \ref{H(2k,n)} below),
and let $L=S$ or $L=Der S$. Then
the corresponding Jordan superalgebras $J(L,\a)$
are as follows:
\begin{enumerate}
\item $L=sl(2m,2n)/\delta_{m,n}\F I_{2m+2n}$, $J(L,\a)=gl(m,n)_+$;
\item $L=osp(m,2r)$, with $m\geq 3$, $(m,r)\neq (4,0)$, $J(L,\a)=(m-3,2r)_+$; 
\item $L=osp(4l,2r)$, with $(l,r)\neq (1,0)$, $J(L,\a)=osp(r,2l)_+$;  
\item $L=p(2k)$ with $k>1$, $J(L,\a)=p(k)_+$;    
\item $L=q(2k)$ with $k> 1$, $J(L,\a)=q(k)_+$;
\item $L=E_7$, $J(L,\a)=E$;    
\item $L=F(4)$, $J(L,\a)=F$;    
\item $L=D(2,1;\alpha)$, $J(L,\a)=D_t$;    
\item $L=Der(sl(2,2)/\F I_4)$, $J(L,\a)=K$.     
\end{enumerate}
\end{remark}

\begin{example}\label{H(2k,n)}\em 
Consider the simple Lie superalgebra $H=H(2k,n)$, where we assume that
$n\geq 3$ if $k>0$ and $n\geq 4$ if $k=0$ (recall that $H(0,n)$ is simple
iff $n\geq 4$). Then
the subalgebra $\a$ of $H$ spanned by elements
$h=\xi_n\xi_{n-1}$, $e=\xi_{n-2}\xi_n$ and $f=\xi_{n-2}\xi_{n-1}$
is a short subalgebra of $H$ whose adjoint action induces on $H$ the
grading $H=H_{-1}\oplus H_0\oplus H_1$
defined by:
\begin{equation}
\begin{array}{c}
\deg(p_i)=\deg(q_i)=0, ~~ i=1, \dots,k;\\
\\
\deg(\xi_j)=0, ~~ j=1, \dots, n-2, ~\deg(\xi_{n-1})=1, 
~\deg(\xi_n)=-1.
\end{array}
\label{degrees}
\end{equation}
Then
$J(H,\a)\simeq {\cal O}(2k,n-2)$ with reversed parity, with the
following product: for $f$, $g$ in $J(H,\a)$,
\begin{equation}
f\circ g=(-1)^{p(f)}\frac{\partial f}{\partial\xi_{n-2}}g+
f\frac{\partial g}{\partial\xi_{n-2}}+(-1)^{p(f)}\xi_{n-2}\{f,g\}
\label{prodH}
\end{equation}
where $\{\cdot, \cdot\}$ is the restriction of bracket (\ref{bracket})
to ${\cal O}(2k, n-2)$.
By Propositions \ref{simplicity}
and \ref{lc}, $J(H,\a)$ is a simple linearly compact Jordan superalgebra.
One easily checks, using (\ref{prodH}), that
the map $$J(H,\a)\longrightarrow JP(2k,n-3),$$  
$$~~~~~~~~~~~~~~~~~~~~~~~~~\xi_{n-2}f+g\longmapsto ~~~~f+\eta g, 
~~~\mbox{for}~f,g\in {\cal O}(2k,n-3),$$ is an isomorphism 
of linearly compact Jordan superalgebras
(see Example \ref{Kdouble}).
\end{example}

\begin{example}\label{K(m,n)}\em Consider the simple Lie superalgebra $L=K(2k+1,n)$,
where we assume that $n\geq 3$. Then
the subalgebra $\a$   spanned by elements
$h=\xi_{n}\xi_{n-1}$, $e=\xi_{n-2}\xi_n$ and $f=\xi_{n-2}\xi_{n-1}$
is a short subalgebra of $L$ whose adjoint action induces on $L$ the
grading $L=L_{-1}\oplus L_0\oplus L_1$
defined by (\ref{degrees}) and the condition:
$\deg(t)=0$. Then
$J(L,\a)\simeq {\cal O}(2k+1,n-2)$ with reversed parity, with the
following product ($f, g \in J(L,\a)$):
\begin{equation}
f\circ g=(-1)^{p(f)}(\frac{\partial f}{\partial\xi_{n-2}}g+
\{\xi_{n-2}f,g\}+\xi_{n-2}(1-E)(f)\frac{\partial g}{\partial t}-
\xi_{n-2}\frac{\partial f}{\partial t}(1-E)(g))
\label{prodK}
\end{equation}
where $\{\cdot, \cdot\}$ is the restriction of bracket (\ref{bracket})
to ${\cal O}(2k+1,n-2)$.
By Propositions \ref{simplicity}
and \ref{lc}, $J(L,\a)$ is a simple linearly compact Jordan superalgebra.
One easily checks, using (\ref{prodK}), that the map
$$J(L,\a)\longrightarrow JP(2k+1,n-3),$$ 
$$~~~~~~~~~~~~~~~~~~~~~~\xi_{n-2}P+Q~~~\longmapsto ~~~~P+\eta Q, 
~~~\mbox{for}~ P, Q\in {\cal O}(2k+1,n-3),$$ 
 is an isomorphism 
of linearly compact Jordan superalgebras (see Example \ref{Kdouble}). 
\end{example}

\begin{example}\label{JS(1,1)}\em Let
$$S'(1,2)=\{X\in W(1,2) ~|~ div(X)=0\},$$
where for $X=P(x,\xi)\frac{\partial}{\partial x}+
\sum_{i=1}^2 Q_i(x,\xi) \frac{\partial}{\partial \xi_i}$,
$div(X)=\frac{\partial P}{\partial x}+
\sum_{i=1}^2 (-1)^{p(Q_i)} \frac{\partial Q_i}{\partial \xi_i}$,
and let $S=S(1,2)=[S'(1,2),S'(1,2)]$.
Then $S'(1,2)=S(1,2)+\F\xi_1\xi_2\frac{\partial}{\partial x}$ and
$S(1,2)$ is simple \cite[Example 4.2]{K2}.
Besides, by \cite[Proposition 6.1]{K3}, $Der S=S\rtimes\mathfrak{a}$
where $\a\simeq sl_2$ acts on $S$ by outer derivations and  is spanned by
$f=\xi_1\xi_2\frac{\partial}{\partial x}+\xi_1\frac{\partial}{\partial\xi_2}$,
$h=\xi_1\frac{\partial}{\partial\xi_1}$ and
$e=-1/2(E+\xi_2\frac{\partial}{\partial\xi_1})$, 
with $E$ defined as follows: for $R=R(x)$,
$E(R\xi_2\frac{\partial}{\partial x}-\frac{\partial R}{\partial x}
\xi_1\xi_2\frac{\partial}{\partial\xi_1})=-R\frac{\partial}{\partial\xi_1}$,
$E(R\xi_1\frac{\partial}{\partial x}+\frac{\partial R}{\partial x}
\xi_1\xi_2\frac{\partial}{\partial\xi_2})=R\frac{\partial}{\partial\xi_2}$,
$E(R\frac{\partial}{\partial x}-
\frac{1}{2}\frac{\partial R}{\partial x}(\xi_1\frac{\partial}{\partial\xi_1}+
\xi_2\frac{\partial}{\partial\xi_2}))=0$,
and $E(R(x)S(0,2))=0$.

Notice that $\a$ is a short
subalgebra of $Der S$ whose action induces on $S$
the $\Z$-grading defined by:
$$\deg x=-\deg\frac{\partial}{\partial x}=0;
~~\deg\xi_1=-\deg\frac{\partial}{\partial \xi_1}=1,
~~\deg\xi_2=-\deg\frac{\partial}{\partial \xi_2}=0.$$ 
Indeed, with this grading,
$S=S_{-1}\oplus S_0\oplus S_1$, where
\begin{itemize}
\item[] $S_{-1}=\langle \frac{\partial}{\partial\xi_1}\rangle
\otimes\F [[x]]\otimes\Lambda(\xi_2)$; 
\item[]
$S_0=(\langle \frac{\partial}{\partial x}, \xi_1\frac{\partial}{\partial\xi_1},
\frac{\partial}{\partial\xi_2}\rangle
\otimes\F [[x]]\otimes\Lambda(\xi_2))\cap S$;
\item[]
$S_1=(\langle \xi_1\frac{\partial}{\partial x}, \xi_1\frac{\partial}{\partial\xi_2}\rangle\otimes\F [[x]]\otimes\Lambda(\xi_2))\cap S$,
\end{itemize}
hence it is immediate to check that $S_i=\{X\in S(1,2)~|~ [h,X]=iX\}$.

By Propositions \ref{simplicity}
and \ref{lc}, $J(S,\a)$ is a simple linearly compact Jordan superalgebra.
An easy computation shows that $J(S,\a)\simeq JS$
(cf.\ Example \ref{new}). 
\end{example}

\begin{example}\label{E(1,6)}\em The Lie superalgebra $L=E(1,6)$
is a simple subalgebra of $K(1,6)$ \cite{K3}. The subalgebra
$\a=\langle \xi_4\xi_6, \xi_6\xi_5, \xi_4\xi_5\rangle$ is a short
subalgebra of $K(1,6)$ contained in $L$, hence it defines on $L$ the
$\Z$-grading  induced by the grading of $K(1,6)$ described in 
Example \ref{K(m,n)} (for $n=6$).
It follows from \cite{KMZ} that $J(L,\a)\simeq JCK$ (cf.\ Example
\ref{JCK}).
\end{example}

\section{Proof of Theorem \ref{main} and 
a Corollary}\label{mainthms}
{\bf Proof of Theorem \ref{main}.} 
By Propositions \ref{simplicity} and  \ref{lc},
 $J$ is a linearly compact unital simple Jordan
superalgebra if and only if $(Lie(J), \a_J)$ is a minimal
pair, with $Lie(J)$ a simple linearly compact Lie superalgebra.
Besides, by Proposition \ref{nonunital}, if $J$ is a non-unital 
simple Jordan superalgebra, then $(Lie(\tilde{J}), \a_{\tilde{J}})$
is a minimal pair with $Lie(\tilde{J})=S\rtimes \a_{\tilde{J}}$,
where $S$ is simple and $\a_{\tilde{J}}$ acts on $S$ by outer derivations.
The result then follows from Theorem \ref{fdshortsa},
Remark \ref{finitedimcases} and Examples \ref{H(2k,n)},
\ref{K(m,n)}, \ref{JS(1,1)} and \ref{E(1,6)}. \hfill$\Box$

\begin{corollary}\label{august} $(a)$ Let $(A, \{\cdot,\cdot\})$ be a unital 
simple linearly compact
generalized Poisson superalgebra. Then the Jordan superalgebra $J(A)$ is
isomorphic to one of the superalgebras $JP(m,n)$.

$(b)$ Any unital simple linearly compact generalized Poisson superalgebra is 
gauge equivalent to one of the generalized Poisson superalgebras $P(m,n)$.

$(c)$ Any unital simple linearly compact Poisson superalgebra is is isomorphic 
to a Poisson superalgebra $P(m,n)$ with $m$ even. 
\end{corollary}
{\bf Proof.}
By Remark \ref{simpleunitalgen}, 
$(A, \{\cdot,\cdot\}_D)$ is simple, hence, by \cite[Theorem 1]{KMCC}, 
its KKM double $J(A)$ is simple.
Furthermore, by construction, $J(A)$ is unital and
admits an involution such that its eigenspace decomposition
$J(A)=J(A)_1\oplus J(A)_{-1}$ has the properties that $J(A)_1$ is
a unital associative subalgebra 
and the right multiplication by its elements induces a $J(A)_1$-module 
structure on $J(A)_{-1}$, isomorphic to the
right regular $J(A)_1$-module with reversed parity. 
Now it is easy to deduce from
Theorem \ref{main} that $J(A)$ is isomorphic to one of the Jordan 
superalgebras $JP(m,n)$, proving $(a)$. 

Since, in the case when
$A$ is purely even, in the above discussion $J(A)_1$
(resp.\ $J(A)_{-1}$) is the even (resp.\ odd) part of $J(A)$,
$(b)$ in the non-super case follows immediately as well.
In the super case we use the description 
of the involutions of $L=Lie(J(A))\simeq H(m,n)$ or $K(m,n)$,
in order to conclude that every involution of $J(A)$
is conjugate to one that changes
signs of indeterminates
(it follows from \cite {CK2} that any involution of $L$ is conjugate to
one of the subgroup $Autgr L \simeq \F^\times(Sp_{2[m/2]}\times O_n)$ of 
$Aut L$). 

Due to the associativity of $J(A)_1$, if
the elements $\eta a_i$ lie in $J(A)_1$, $i=1,2,3$, then
\begin{equation}
\{a_1,a_2\}a_3=a_1\{a_2, a_3\}.
\label{associativity}
\end{equation}
We check, using condition (\ref{associativity}), that $J(A)_1=A$ and
$J(A)_{-1}=\eta A$ for all $A=P(m,n)$, except for $A=P(0,1)$,
in which case
there is one other involution of $J(A)$
with required properties, but it is conjugate to the standard one
by the permutation of $\xi_1$ and $\eta$. This proves $(b)$.
Claim $(c)$ follows easily from $(b)$. 
\hfill $\Box$

\medskip

A classification of local transitive Lie algebras on a 
finite-dimensional
manifold in terms, analogous to those of Corollary \ref{august}$(b)$,
was obtained in \cite{Ki}. The finite-dimensional 
case of Corollary \ref{august}$(c)$ was obtained in \cite{Ka2}
(the argument there is not quite correct).

\begin{remark}\em Gauge equivalence classes of generalized Poisson brackets 
on a 
unital associative commutative algebra $A$ are in a canonical bijective
correspondence with Jordan superalgebras $J$, such that $J_{\bar{0}}\simeq A$,
and the right multiplication by elements of $J_{\bar{0}}$ defines an 
$A$-module 
structure on $J_{\bar{1}}$, isomorphic to the right regular $A$-module.
\end{remark}


\setcounter{section}{1}
\renewcommand{\thesection}{\Alph{section}}
\setcounter{theorem}{0}
\section*{Appendix A. All known infinite-dimensional simple Jordan 
superalgebras over $\F$.}
\begin{examples}\em All known examples of infinite-dimensional 
simple Jordan superalgebras
over an algebraically closed field $\F$ of characteristic
zero, are of the following type:
\begin{enumerate}
\item $A^+$, where $A$ is a simple associative superalgebra,
with Jordan product:
$$a\circ b=\frac{1}{2}(ab+(-1)^{p(a)p(b)}ba).$$
\item $BCP(A,\omega)$: the fixed point subalgebra of $A^+$ under a super
anti-involution $\omega$.
\item $D(V, (\cdot,\cdot))$: the Jordan superalgebra of a non degenerate
super-symmetric bilinear form.
\item $Q(A)=A\oplus\bar{A}$, where $A$ is a simple associative algebra,
$\bar{A}$ a copy of $A$ with odd parity, with Jordan product:
$$(a\oplus b)\circ (a_1\oplus b_1)=(aa_1+bb_1, ab_1+ba_1).$$
\item $J(P)$: the KKM double of a simple generalized Poisson superalgebra $P$.
\item $JS(A,D)$, where $A$ is a commutative superalgebra with an odd 
derivation $D$, such that $A$ contains no non-trivial $D$-invariant ideals,
with reversed parity and Jordan product:
$$a\circ b  =aD(b)+(-1)^{p(a)}D(a)b.$$
(The example $JS(\Lambda(1),d/d\xi)$ shows that these Jordan superalgebras
are not always simple.)
\item $JCK(A,D)$, where $A$ is a commutative associative algebra with an
even derivation $D$, such that $A$ contains no non-trivial $D$-invariant
 ideals (see Example \ref{JCK}).
\end{enumerate}
\end{examples}

\begin{question}\em Are there any other examples?
\end{question}

\begin{question}\em It is shown in \cite{RZ} that if char $\F>3$, all
unital finite-dimensional simple Jordan superalgebras over $\F$ are either
of types 1.-5., 7., or they are isomorphic to $D_t$, $E$, $F$. Is it true
that all the non-unital ones are either of type 6. or isomorphic to $K$?
\end{question}

The following remarks explain why there are no infinite-dimensional
linearly compact simple Jordan superalgebras of types 1-4.

\begin{remarks}\em (a) Any linearly compact associative superalgebra
contains an ideal of finite codimension. Hence the simple
ones are finite-dimensional.

(b) A non-degenerate bilinear form exists on a linearly compact space $V$
(resp.\ $Hom(V,V)$ is linearly compact) if and only if $\dim V<\infty$.
\end{remarks}

\section*{Appendix B. The finite-dimensional Jordan 
superalgebras
$\boldsymbol{D_t}$, $\boldsymbol{K}$ and $\boldsymbol{F}$.}
In this section we recall the construction of the exceptional
 finite-dimensional Jordan superalgebras
$D_t$, $K$ and $F$.

\medskip

\noindent
{\bf The superalgebra $\boldsymbol{D_t}$ ($\boldsymbol{t\in\F, t\neq 0}$) .}
This is a parametric series 
of unital four-dimensional Jordan superalgebras defined as follows:
$(D_t)_{\bar{0}}=\F e_1+\F e_2$, $(D_t)_{\bar{1}}=\F \xi +\F \eta$,
with multiplication table:
$$e_i^2=e_i,~~e_1\circ e_2=0,~~e_i\circ \xi=\frac{1}{2}\xi,
~~e_i\circ \eta=\frac{1}{2}\eta, ~~\xi\circ \eta=e_1+te_2.$$
$D_t$ is simple if and only if $t\neq 0$.

\medskip

\noindent
{\bf The superalgebra $\boldsymbol{K}$.} This is the three-dimensional simple
non-unital Jordan superalgebra 
$K=\langle a\rangle\oplus\langle\xi_1,\xi_2\rangle$,
with multiplication table:
$$a^2=a, ~~a\circ\xi_i=\frac{1}{2}\xi_i,~~\xi_1\circ\xi_2=a.$$

\medskip

\noindent
{\bf The superalgebra $\boldsymbol{F}$.} 
In \cite{K1}, \cite{HK} it is defined by a multiplication table. More
recently, I.\ Shestakov proposed the following beautiful construction
(see \cite{RZ}).
Consider on $K$ the supersymmetric
bilinear form $(\cdot, \cdot)$ defined by:
$$(a,a)=\frac{1}{2}, ~~(\xi_1, \xi_2)=1, ~~(a, \xi_i)=0.$$
Then $F=\F 1+(K\otimes K)$, with the unit element 1 and the product:
$$(a\otimes b)\circ (c\otimes d)=(-1)^{p(b)p(c)}(ac\otimes bd-
\frac{3}{4}(a,c)(b,d)1),$$
for $a, b, c, d\in K$.

$$$$ 

\end{document}